\documentclass[A4j,11pt]{article}
\usepackage{graphics}
\usepackage{amsmath}
\usepackage{amssymb}
\usepackage{latexsym}
\usepackage{amsfonts}
\begin{document}

\title{{\bf Mean curvature flow of\\
certain kind of isoparametric foliations\\
on non-compact symmetric spaces}}
\author{{\bf Naoyuki Koike}}
%
%
\date{}
%
\maketitle
\begin{abstract}
In this paper, we investigate the mean curvature flows starting from all non-minimal leaves of the isoparametric 
foliation given by a certain kind of solvable group action on a symmetric space of non-compact type.  
We prove that the mean curvature flow starting from each non-minimal leaf of the foliation 
exists in infinite time, if the foliation admits no minimal leaf, then the flow asymptotes 
the self-similar flow starting from another leaf, and if the foliation admits a minimal leaf 
(in this case, it is shown that there exists the only one minimal leaf), 
then the flow converges to the minimal leaf of the foliation in $C^{\infty}$-topology.  
\end{abstract}

\vspace{0.5truecm}

\section{Introduction}
In [K2], we proved that the mean curvature flow starting from any non-minimal compact isoparametric 
(equivalently, equifocal) submanifold in a symmetric space of compact type collapses to one of its focal 
submanifolds in finite time.  
Here we note that parallel submanifolds and focal ones of the isoparametric submanifold give 
an isoparametric foliation consisting of compact leaves on the symmetric space, where 
an {\it isoparametric foliation} means a singular Riemannian foliation satisfying the following conditions:

\vspace{0.2truecm}

(i) The mean curvature form is basic, 

(ii) The regular leaves are submanifolds with section.  

\vspace{0.2truecm}

\noindent
A singular Riemannian foliation satisfying only the first condition is called a {\it generalized isoparametric 
foliation}.  Recently, M. M. Alexandrino and M. Radeschi ([AR]) investigated the mean curvature flow 
starting from a regular leaf of a generalized isoparametric foliation consisting of compact leaves on a compact 
Riemannian manifold.  In particular, they ([AR]) generalized our result to the mean curvature flow 
starting from a regular leaf of the foliation in the case where the foliation is isoparametric and 
the ambient space curves non-negatively.  
On the other hand, we ([K3]) proved that the mean curvature flow starting from a certain kind of non-minimal 
(not necessarily compact) isoparametric submanifold in a symmetric space of non-compact type 
(which curves non-positively) collapses to one of its focal submanifolds in finite time.  
Here we note that the isoparametric foliation associated with this isoparametric submanifold 
consists of curvature-adapted leaves.  See the next paragraph about the definition of the curvature-adaptedness.  

In this paper, we study the mean curvature flow starting from leaves of the isoparametiric foliation given by 
the action of a certain kind of solvable subgroup (see Examples 1 and 2) of the (full) isometry group of 
a symmetric space of non-compact type.  
Here we note that this isoparametric foliation consists of (not necessarily curvature-adapted) 
non-compact regular leaves.  
We shall explain the solvable group action which we treat in this paper.  
Let $G/K$ be a symmetric space of non-compact type, $\mathfrak g=\mathfrak k+\mathfrak p$ 
($\mathfrak k:={\rm Lie}\,K$) be the Cartan decomposition associated with the symmetric pair $(G,K)$, 
$\mathfrak a$ be the maximal abelian subspace of $\mathfrak p$, $\widetilde{\mathfrak a}$ be the Cartan 
subalgebra of $\mathfrak g$ containing $\mathfrak a$ and $\mathfrak g=\mathfrak k+\mathfrak a+\mathfrak n$ be 
the Iwasawa's decomposition.  Let $A,\,\widetilde A$ and $N$ be the connected Lie subgroups of $G$ having 
$\mathfrak a,\,\widetilde{\mathfrak a}$ and $\mathfrak n$ as their Lie algebras, respectively.  
Let $\pi:G\to G/K$ be the natural projection.  

\vspace{0.3truecm}

\noindent
{\bf Given metric.} In this paper, we give $G/K$ the $G$-invariant metric induced from the restriction 
$B\vert_{\mathfrak p\times\mathfrak p}$ of the Killing form $B$ of $\mathfrak g$ to 
$\mathfrak p\times\mathfrak p$.  

\vspace{0.3truecm}

\noindent
The symmetric space $G/K$ is identified with the solvable group $AN$ with a left-invariant metric through 
$\pi\vert_{AN}$.  Fix a lexicographic ordering of $\mathfrak a$.  
Let $\mathfrak g=\mathfrak g_0+\sum\limits_{\lambda\in\triangle}
\mathfrak g_{\lambda}$, $\mathfrak p=\mathfrak a
+\sum\limits_{\lambda\in\triangle_+}\mathfrak p_{\lambda}$ and 
$\mathfrak k=\mathfrak k_0+\sum\limits_{\lambda\in\triangle_+}
\mathfrak k_{\lambda}$ be the root space decompositions of 
$\mathfrak g,\,\mathfrak p$ and $\mathfrak k$ with respect to $\mathfrak a$, 
where we note that 
$$\begin{array}{l}
\displaystyle{\mathfrak g_{\lambda}=\{X\in\mathfrak g\,\vert\,{\rm ad}(a)X
=\lambda(a)X\,\,{\rm for}\,\,{\rm all}\,\,a\in\mathfrak a\}\,\,\,\,
(\lambda\in\triangle),}\\
\displaystyle{\mathfrak p_{\lambda}=\{X\in\mathfrak p\,\vert\,{\rm ad}(a)^2X
=\lambda(a)^2X\,\,{\rm for}\,\,{\rm all}\,\,a\in\mathfrak a\}\,\,\,\,
(\lambda\in\triangle_+),}\\
\displaystyle{\mathfrak k_{\lambda}=\{X\in\mathfrak k\,\vert\,{\rm ad}(a)^2X
=\lambda(a)^2X\,\,{\rm for}\,\,{\rm all}\,\,a\in\mathfrak a\}\,\,\,\,
(\lambda\in\triangle_+\cup\{0\}).}
\end{array}$$
Note that $\mathfrak n=\sum\limits_{\lambda\in\triangle_+}\mathfrak g_{\lambda}$.  
Let $G=KAN$ be the Iwasawa decomposition of $G$.  
Now we shall give examples of a solvable group contained in $AN$ whose 
action on $G/K(=AN)$ is (complex) hyperpolar.  
Since $G/K$ is of non-compact type, 
$\pi$ gives a diffeomorphism of $AN$ onto $G/K$.  
Denote by $\langle\,\,,\,\,\rangle$ the left-invariant metric of $AN$ induced 
from the metric of $G/K$ by $\pi\vert_{AN}$.  Also, denote by $\langle\,\,,\,\,
\rangle^G$ the bi-invariant metric of $G$ induced from the Killing form $B$.  Note that 
$\langle\,\,,\,\,\rangle\not=\iota^{\ast}\langle\,\,,\,\,\rangle^G$, where 
$\iota$ is the inclusion map of $AN$ into $G$.  
Denote by ${\rm Exp}$ the exponential map of the Riemannian manifold $AN(=G/K)$ at $e$ and 
by $\exp_G$ the exponential map of the Lie group $G$.  
Let ${\it l}$ be a $r$-dimensional subspace of $\mathfrak a+\mathfrak n$ and 
set $\mathfrak s:=(\mathfrak a+\mathfrak n)\ominus{\it l}$, where 
$(\mathfrak a+\mathfrak n)\ominus {\it l}$ denotes the orthogonal complement 
of ${\it l}$ in $\mathfrak a+\mathfrak n$ with respect to 
$\langle\,\,,\,\,\rangle_e$ ($e\,:\,$ is the identity element of $G$).  
According to the result in [K1], if $\mathfrak s$ is a subalgebra of $\mathfrak a+\mathfrak n$ and 
${\it l}_{\mathfrak p}:={\rm pr}_{\mathfrak p}({\it l})$ (${\rm pr}_{\mathfrak p}\,:\,$ the orthogonal projection 
of $\mathfrak g$ onto $\mathfrak p$) is abelian, then the $S$-action ($S:=\exp_G(\mathfrak s)$) 
gives an isoparametric foliation without singular leaf.  
We ([K1]) gave examples of such a subalgebra $\mathfrak s$ of $\mathfrak a+\mathfrak n$.  

\vspace{0.5truecm}

\noindent
{\it Example 1.} Let $\mathfrak b$ be a $r(\geq1)$-dimensional subspace of 
$\mathfrak a$ and $\mathfrak s_{\mathfrak b}:=(\mathfrak a+\mathfrak n)\ominus\mathfrak b$.  
It is clear that $\mathfrak b_{\mathfrak p}(=\mathfrak b)$ is abelian and that $\mathfrak s_{\mathfrak b}$ 
is a subalgebra of $\mathfrak a+\mathfrak n$.  

\vspace{0.5truecm}

\noindent
{\it Example 2.} Let $\{\lambda_1,\cdots,\lambda_k\}$ be a subset of a simple 
root system $\Pi$ of $\triangle$ such that $H_{\lambda_1},\cdots,
H_{\lambda_k}$ are mutually orthogonal, $\mathfrak b$ be a subspace of 
$\mathfrak a\ominus{\rm Span}\{H_{\lambda_1},\cdots,H_{\lambda_k}\}$ 
(where $\mathfrak b$ may be $\{0\}$) and ${\it l}_i$ ($i=1,\cdots,k$) be 
a one-dimensional subspace of ${\bf R}H_{\lambda_i}+\mathfrak g_{\lambda_i}$ with 
${\it l}_i\not={\bf R}H_{\lambda_i}$, where $H_{\lambda_i}$ is the 
element of $\mathfrak a$ defined by $\langle H_{\lambda_i},\cdot\rangle=
\lambda_i(\cdot)$ and ${\bf R}H_{\lambda_i}$ is the subspace of $\mathfrak a$ 
spanned by $H_{\lambda_i}$.  Set ${\it l}:=\mathfrak b+\sum\limits_{i=1}^k
{\it l}_i$.  
Then, it is shown that ${\it l}_{\mathfrak p}$ is abelian and that 
$\mathfrak s_{\mathfrak b,{\it l}_1,\cdots,{\it l}_k}:=(\mathfrak a+\mathfrak n)\ominus{\it l}$ 
is a subalgebra of $\mathfrak a+\mathfrak n$.  

\vspace{0.5truecm}

\noindent
In Example 2, a unit vector of ${\it l}_i$ is described as 
$\displaystyle{\frac{1}{\cosh(\vert\vert\lambda_i\vert\vert t_i)}\xi^i-\frac{1}{\vert\vert\lambda_i\vert\vert}
\tanh(\vert\vert\lambda_i\vert\vert t_i)H_{\lambda_i}}$ for a unit vector $\xi^i$ of 
$\mathfrak g_{\lambda_i}$ and some $t_i\in{\mathbb R}$, where 
$\vert\vert\lambda_i\vert\vert:=\vert\vert H_{\lambda_i}\vert\vert$.  
Then we denote ${\it l}_i$ by ${\it l}_{\xi^i,t_i}$ if necessary and 
set $\displaystyle{\xi^i_{t_i}:=\frac{1}{\cosh(\vert\vert\lambda_i\vert\vert t_i)}\xi^i
-\frac{1}{\vert\vert\lambda_i\vert\vert}\tanh(\vert\vert\lambda_i\vert\vert t_i)H_{\lambda_i}}$.  
Set $S_{\mathfrak b}:=\exp_G(\mathfrak s_{\mathfrak b})$ and 
$S_{\mathfrak b,,{\it l}_1,\cdots,{\it l}_k}:=\exp_G(\mathfrak s_{\mathfrak b,{\it l}_1,\cdots,{\it l}_k})$.  
Denote by ${\mathfrak F}_{\mathfrak b}$ and ${\mathfrak F}_{\mathfrak b,{\it l}_1,\cdots,{\it l}_k}$ 
the isoparametric foliations given by the $S_{\mathfrak b}$-action and 
the $S_{\mathfrak b,{\it l}_1,\cdots,{\it l}_k}$-one, respectively.  
A submanifold in a Riemannian manifold is said to be {\it curvature-adapted} if, for each normal vector $v$ of 
the submanifold, the normal Jacobi operator $R(v):=R(\cdot,v)v$ preserves the tangent space 
of the submanifold invariantly and the restriction of $R(v)$ to the tangent space commutes with the shape 
operator $A_v$, where $R$ is the curvature tensor of the ambient Riemannian manifold.  
According to the results in [K1], the following facts hold for isoparametric foliations 
${\mathfrak F}_{\mathfrak b}$ and ${\mathfrak F}_{\mathfrak b,{\it l}_1,\cdots,{\it l}_k}$:

\vspace{0.2truecm}

{\rm (i)} All leaves of ${\mathfrak F}_{\mathfrak b}$ are curvature-adapted.  

{\rm (ii)} Let $\lambda_1,\cdots,\lambda_k\,(\in\triangle_+)$ be as in Example 2.  
If the root system $\triangle$ of $G/K$ is non-reduced and 
$2\lambda_{i_0}\in\triangle_+$ for some $i_0\in\{1,\cdots,k\}$, then all leaves of 
${\mathfrak F}_{\mathfrak b,{\it l}_1,\cdots,{\it l}_k}$ are not curvature-adapted.  

{\rm (iii)} If $\mathfrak b\not=\{0\}$, then ${\mathfrak F}_{\mathfrak b,{\it l}_1,\cdots,{\it l}_k}$ 
admits no minimal leaf.  
On the other hand, if $\mathfrak b=\{0\}$, then this action admits the only minimal leaf.  

{\rm (iv)} Let ${\it l}_1,\cdots,{\it l}_k$ be as in Example 2 and 
$\overline{\it l}_i$ ($i=1,\cdots,k$) be the orthogonal projection of ${\it l}_i$ onto 
$\mathfrak g_{\lambda_i}$.  
Then $\mathfrak F_{\mathfrak b,\overline{\it l}_1,\cdots,\overline{\it l}_k}$ is congruent to 
${\mathfrak F}_{\mathfrak b,{\it l}_1,\cdots,{\it l}_k}$.  
In more detail, we have 
$$L_{b\cdot\gamma_{\xi^1}(t_1)\cdot\,\cdots\,\cdot\gamma_{\xi^k}(t_k)}
(S_{\mathfrak b,{\it l}_1,\cdots,{\it l}_k}\cdot e)
=S_{\mathfrak b,\overline{\it l}_1,\cdots,\overline{\it l}_k}\cdot
(b\cdot\gamma_{\xi^1}(t_1)\cdot\,\cdots\,\cdot\gamma_{\xi^k}(t_k)),$$
where $\gamma_{\xi^i}$ ($i=1,\cdots,k$) is the geodesic in $AN(=G/K)$ with $\gamma_{\xi^i}'(0)=\xi^i$, 
$b$ is an element of $\exp(\mathfrak b)$ and 
$L_{b\cdot\gamma_{\xi^1}(t_1)\cdot\,\cdots\,\cdot\gamma_{\xi^k}(t_k)}$ is the left translation by 
$b\cdot\gamma_{\xi^1}(t_1)\cdot\,\cdots\,\cdot\gamma_{\xi^k}(t_k)$.  
For example, in case of $k=1$ and $b=e$, the positional relation among the leaves of these foliations is 
as in Figure 1.  

\vspace{0.5truecm}

\centerline{
\unitlength 0.1in
\begin{picture}( 90.9000, 12.9000)(-54.9000,-21.6000)
%
\special{pn 8}%
\special{pa 1600 2000}%
\special{pa 3200 2000}%
\special{fp}%
%
\special{pn 8}%
\special{pa 1600 2000}%
\special{pa 2000 1600}%
\special{fp}%
%
\special{pn 8}%
\special{ar 2430 1800 230 980  3.1415927 4.0031485}%
%
\special{pn 8}%
\special{pa 2200 1800}%
\special{pa 2186 1772}%
\special{pa 2172 1744}%
\special{pa 2158 1714}%
\special{pa 2144 1686}%
\special{pa 2132 1656}%
\special{pa 2118 1628}%
\special{pa 2106 1598}%
\special{pa 2094 1568}%
\special{pa 2080 1538}%
\special{pa 2070 1508}%
\special{pa 2056 1480}%
\special{pa 2046 1450}%
\special{pa 2034 1420}%
\special{pa 2022 1390}%
\special{pa 2012 1360}%
\special{pa 2002 1328}%
\special{pa 1992 1298}%
\special{pa 1982 1268}%
\special{pa 1974 1238}%
\special{pa 1966 1206}%
\special{pa 1956 1176}%
\special{pa 1950 1144}%
\special{pa 1942 1114}%
\special{pa 1940 1100}%
\special{sp}%
%
\special{pn 20}%
\special{sh 1}%
\special{ar 2200 1800 10 10 0  6.28318530717959E+0000}%
\special{sh 1}%
\special{ar 2200 1800 10 10 0  6.28318530717959E+0000}%
%
\special{pn 8}%
\special{pa 1990 2000}%
\special{pa 2390 1600}%
\special{fp}%
%
\special{pn 8}%
\special{pa 3200 2000}%
\special{pa 3600 1600}%
\special{fp}%
%
\special{pn 8}%
\special{pa 1660 2150}%
\special{pa 3016 1544}%
\special{fp}%
%
\special{pn 8}%
\special{pa 1650 2160}%
\special{pa 1736 2038}%
\special{fp}%
%
\special{pn 8}%
\special{pa 2010 1600}%
\special{pa 2070 1600}%
\special{fp}%
%
\special{pn 8}%
\special{pa 2130 1600}%
\special{pa 2180 1600}%
\special{fp}%
%
\special{pn 8}%
\special{pa 2240 1600}%
\special{pa 2790 1600}%
\special{fp}%
%
\special{pn 8}%
\special{pa 2940 1600}%
\special{pa 3600 1600}%
\special{fp}%
%
\special{pn 8}%
\special{pa 2200 1790}%
\special{pa 2390 1790}%
\special{fp}%
%
\special{pn 8}%
\special{pa 2510 1790}%
\special{pa 3410 1790}%
\special{fp}%
%
\special{pn 8}%
\special{pa 2200 1790}%
\special{pa 3160 1360}%
\special{fp}%
%
\special{pn 8}%
\special{ar 3530 1790 230 980  3.1415927 4.0031485}%
%
\special{pn 20}%
\special{sh 1}%
\special{ar 3300 1790 10 10 0  6.28318530717959E+0000}%
\special{sh 1}%
\special{ar 3300 1790 10 10 0  6.28318530717959E+0000}%
\put(20.5000,-10.9000){\makebox(0,0)[rb]{$S_{\mathfrak b,{\it l}_1}\cdot e$}}%
\put(21.0000,-10.5000){\makebox(0,0)[lb]{$S_{\mathfrak b,\overline{\it l}_1}\cdot e$}}%
\put(51.3000,-10.4000){\makebox(0,0)[rb]{$S_{\mathfrak b,\overline{\it l}_1}\cdot\gamma_{\xi^1}(t_1)=L_{\gamma_{\xi^1}(t_1)}(S_{\mathfrak b,{\it l}_1}\cdot e)$}}%
%
\special{pn 8}%
\special{pa 1590 1650}%
\special{pa 1750 1950}%
\special{dt 0.045}%
\special{sh 1}%
\special{pa 1750 1950}%
\special{pa 1736 1882}%
\special{pa 1726 1904}%
\special{pa 1702 1902}%
\special{pa 1750 1950}%
\special{fp}%
%
\special{pn 8}%
\special{pa 1470 2010}%
\special{pa 1790 2060}%
\special{dt 0.045}%
\special{sh 1}%
\special{pa 1790 2060}%
\special{pa 1728 2030}%
\special{pa 1738 2052}%
\special{pa 1722 2070}%
\special{pa 1790 2060}%
\special{fp}%
\put(6.7000,-20.9000){\makebox(0,0)[lb]{${\rm Exp}(\mathfrak b+{\it l}_1)$}}%
\put(12.0000,-16.2000){\makebox(0,0)[lb]{${\rm Exp}(\mathfrak b+\overline{\it l}_1)$}}%
\put(29.1000,-23.0000){\makebox(0,0)[lb]{$\gamma_{\xi^1}$}}%
\put(35.4000,-14.0000){\makebox(0,0)[lb]{$\gamma_{\xi^1}(t_1)$}}%
%
\special{pn 8}%
\special{pa 2980 2090}%
\special{pa 2920 1790}%
\special{dt 0.045}%
\special{sh 1}%
\special{pa 2920 1790}%
\special{pa 2914 1860}%
\special{pa 2930 1842}%
\special{pa 2954 1852}%
\special{pa 2920 1790}%
\special{fp}%
\put(27.1000,-12.9000){\makebox(0,0)[rb]{$\gamma_{\xi^1_{t_1}}$}}%
\put(21.3000,-18.4000){\makebox(0,0)[rb]{$e$}}%
%
\special{pn 8}%
\special{pa 2670 1320}%
\special{pa 2820 1510}%
\special{dt 0.045}%
\special{sh 1}%
\special{pa 2820 1510}%
\special{pa 2794 1446}%
\special{pa 2788 1468}%
\special{pa 2764 1470}%
\special{pa 2820 1510}%
\special{fp}%
%
\special{pn 8}%
\special{pa 2390 1600}%
\special{pa 3280 1210}%
\special{fp}%
%
\special{pn 8}%
\special{pa 3010 1550}%
\special{pa 3280 1220}%
\special{fp}%
%
\special{pn 8}%
\special{pa 2380 1600}%
\special{pa 2224 1660}%
\special{fp}%
%
\special{pn 8}%
\special{pa 2180 1660}%
\special{pa 2160 1670}%
\special{fp}%
%
\special{pn 8}%
\special{pa 2120 1700}%
\special{pa 1906 1780}%
\special{fp}%
%
\special{pn 8}%
\special{pa 1790 1980}%
\special{pa 1910 1780}%
\special{fp}%
%
\special{pn 8}%
\special{pa 3560 1460}%
\special{pa 3320 1770}%
\special{dt 0.045}%
\special{sh 1}%
\special{pa 3320 1770}%
\special{pa 3378 1730}%
\special{pa 3354 1728}%
\special{pa 3346 1706}%
\special{pa 3320 1770}%
\special{fp}%
\end{picture}%
\hspace{19.5truecm}}

\vspace{0.5truecm}

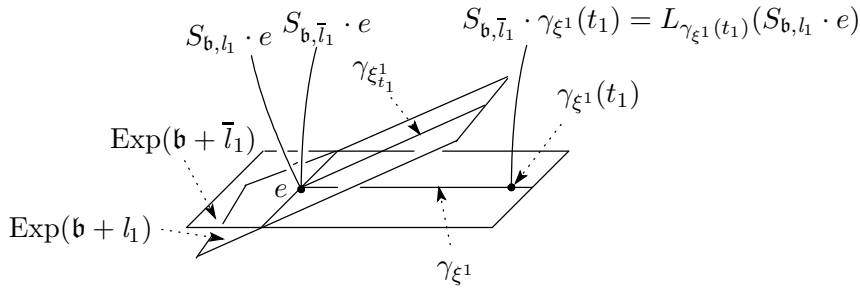
\centerline{{\bf Figure 1.}}

\vspace{0.5truecm}

\noindent
According to the above facts (i) and (ii), the leaves of ${\mathfrak F}_{\mathfrak b,{\it l}_1,\cdots,{\it l}_k}$ 
give examples of interesting isoparametric submanifolds in $G/K$.  

In this paper, we shall prove the following facts for the mean curvature flows starting from the non-minimal 
leaves of ${\mathfrak F}_{\mathfrak b,\overline{\it l}_1,\cdots,\overline{\it l}_k}$.  

\vspace{0.3truecm}

\noindent
{\bf Theorem A.} {\sl Assume that $\mathfrak b\not=\{0\}$.  
Let $M$ be any leaf of ${\mathfrak F}_{\mathfrak b,\overline{\it l}_1,\cdots,\overline{\it l}_k}$.  
and $M_t$ ($0\leq t<T$) be the mean curvature flow starting from $M$.  
Then the following statements ${\rm (i)}-{\rm(iii)}$ hold.  

{\rm (i)} $T=\infty$ holds.  

{\rm (ii)} If $M$ passes through $\exp(\mathfrak b)$, then the mean curvature flow $M_t$ is self-similar.  

{\rm (iii)} If $M$ does not pass through $\exp(\mathfrak b)$, then the mean curvature flow $M_t$ asymptotes 
the mean curvature flow starting from the leaf of 
${\mathfrak F}_{\mathfrak b,\overline{\it l}_1,\cdots,\overline{\it l}_k}$ passing through a point of 
$\exp(\mathfrak b)$.}

\vspace{0.3truecm}

\noindent
{\it Remark 1.1.} The mean curvature flow starting from any leaf of ${\mathfrak F}_{\mathfrak b}$ 
is self-similar.  

\vspace{0.3truecm}

\centerline{
\unitlength 0.1in
\begin{picture}( 41.1000, 17.1000)( -9.1000,-20.8000)
%
\special{pn 8}%
\special{pa 1600 1000}%
\special{pa 1210 1600}%
\special{pa 2800 1600}%
\special{pa 3200 1000}%
\special{pa 3200 1000}%
\special{pa 1600 1000}%
\special{fp}%
%
\special{pn 8}%
\special{pa 2020 1600}%
\special{pa 2380 1000}%
\special{fp}%
%
\special{pn 8}%
\special{pa 2630 1480}%
\special{pa 2610 1456}%
\special{pa 2590 1432}%
\special{pa 2570 1406}%
\special{pa 2552 1380}%
\special{pa 2534 1354}%
\special{pa 2518 1326}%
\special{pa 2502 1298}%
\special{pa 2486 1270}%
\special{pa 2472 1242}%
\special{pa 2460 1212}%
\special{pa 2448 1182}%
\special{pa 2436 1152}%
\special{pa 2424 1122}%
\special{pa 2416 1092}%
\special{pa 2406 1062}%
\special{pa 2400 1030}%
\special{pa 2394 998}%
\special{pa 2394 996}%
\special{sp -0.045}%
%
\special{pn 8}%
\special{pa 2330 1000}%
\special{pa 2316 1030}%
\special{pa 2302 1058}%
\special{pa 2284 1084}%
\special{pa 2264 1108}%
\special{pa 2242 1132}%
\special{pa 2220 1154}%
\special{pa 2194 1176}%
\special{pa 2168 1194}%
\special{pa 2142 1212}%
\special{pa 2114 1226}%
\special{pa 2086 1242}%
\special{pa 2056 1254}%
\special{pa 2026 1264}%
\special{pa 1996 1276}%
\special{pa 1964 1284}%
\special{pa 1934 1292}%
\special{pa 1902 1298}%
\special{pa 1870 1304}%
\special{pa 1840 1308}%
\special{pa 1808 1310}%
\special{pa 1776 1312}%
\special{pa 1772 1312}%
\special{sp -0.045}%
%
\special{pn 8}%
\special{ar 2290 1520 220 890  3.1415927 4.1226726}%
%
\special{pn 8}%
\special{ar 2000 1310 220 890  3.1415927 4.1226726}%
%
\special{pn 8}%
\special{ar 2420 1480 220 890  5.3021054 6.2831853}%
%
\special{pn 13}%
\special{ar 1790 940 540 360  1.0841402 1.5794541}%
%
\special{pn 13}%
\special{pa 2460 1200}%
\special{pa 2440 1150}%
\special{fp}%
\special{sh 1}%
\special{pa 2440 1150}%
\special{pa 2446 1220}%
\special{pa 2460 1200}%
\special{pa 2484 1204}%
\special{pa 2440 1150}%
\special{fp}%
%
\special{pn 13}%
\special{ar 3600 940 1190 930  2.5416669 2.8662134}%
%
\special{pn 13}%
\special{ar 1780 830 620 470  0.7645195 0.9918982}%
%
\special{pn 13}%
\special{pa 2220 1160}%
\special{pa 2260 1120}%
\special{fp}%
\special{sh 1}%
\special{pa 2260 1120}%
\special{pa 2200 1154}%
\special{pa 2222 1158}%
\special{pa 2228 1182}%
\special{pa 2260 1120}%
\special{fp}%
\put(26.6000,-12.1000){\makebox(0,0)[lb]{${\rm Exp}({\it l})$}}%
%
\special{pn 8}%
\special{pa 1300 1160}%
\special{pa 2040 1550}%
\special{dt 0.045}%
\special{sh 1}%
\special{pa 2040 1550}%
\special{pa 1990 1502}%
\special{pa 1994 1526}%
\special{pa 1972 1538}%
\special{pa 2040 1550}%
\special{fp}%
\put(12.5000,-11.5000){\makebox(0,0)[rb]{${\rm Exp}(\mathfrak b)$}}%
\put(18.4000,-5.4000){\makebox(0,0)[lb]{$M^1$}}%
\put(25.1000,-7.2000){\makebox(0,0)[lb]{$M^3$}}%
\put(21.4000,-7.5000){\makebox(0,0)[lb]{$M^2$}}%
\put(4.1000,-17.0000){\makebox(0,0)[lt]{{\small The mean curvature flows starting from leaves $M^1$ and $M^3$}}}%
%
\special{pn 13}%
\special{pa 2070 1520}%
\special{pa 2250 1210}%
\special{fp}%
\special{sh 1}%
\special{pa 2250 1210}%
\special{pa 2200 1258}%
\special{pa 2224 1256}%
\special{pa 2234 1278}%
\special{pa 2250 1210}%
\special{fp}%
\put(4.1000,-18.9000){\makebox(0,0)[lt]{{\small of ${\mathfrak F}_{\mathfrak b,\bar{\it l}_1,\cdots,\bar{\it l}_k}$ ($\mathfrak b\not=\{0\}$) asymptotes the mean curvature flow}}}%
\put(4.1000,-20.8000){\makebox(0,0)[lt]{{\small (which is self-similar) starting from a leaf $M^2$ of ${\mathfrak F}_{\mathfrak b,\bar{\it l}_1,\cdots,\bar{\it l}_k}$.}}}%
\end{picture}%
\hspace{5.5truecm}}

\vspace{0.5truecm}

\centerline{{\bf Figure 2.}}

\vspace{0.3truecm}

Also, in case of $\mathfrak b=\{0\}$, we obtain the following fact.  

\vspace{0.3truecm}

\noindent
{\bf Theorem B.} {\sl 
Let $M$ be a leaf of ${\mathfrak F}_{\{0\},\overline{\it l}_1,\cdots,\overline{\it l}_k}$-action other than 
$S_{\{0\},\overline{\it l}_1,\cdots,\overline{\it l}_k}\cdot e$ and $M_t$ ($0\leq t<T$) be the mean curvature 
flow starting from $M$.  Then the following statements ${\rm (i)}-{\rm(ii)}$ hold.  

{\rm (i)} $T=\infty$ holds.  

{\rm (ii)} $M_t$ convergres to the only minimal leaf 
$S_{\{0\},\overline{\it l}_1,\cdots,\overline{\it l}_k}\cdot e$ (in $C^{\infty}$-topology) as $t\to\infty$.  
}

\vspace{0.2truecm}

\centerline{
\unitlength 0.1in
\begin{picture}( 28.8000, 16.6000)(  3.2000,-19.1000)
%
\special{pn 8}%
\special{pa 1600 1000}%
\special{pa 1210 1600}%
\special{pa 2800 1600}%
\special{pa 3200 1000}%
\special{pa 3200 1000}%
\special{pa 1600 1000}%
\special{fp}%
%
\special{pn 8}%
\special{ar 2210 1480 220 890  3.1415927 4.1226726}%
%
\special{pn 8}%
\special{ar 2030 1170 220 890  3.1415927 4.1226726}%
%
\special{pn 8}%
\special{ar 2380 1230 220 890  5.3021054 6.2831853}%
\put(23.8000,-15.5000){\makebox(0,0)[lb]{${\rm Exp}({\it l})$}}%
\put(18.7000,-4.2000){\makebox(0,0)[lb]{$M^1$}}%
\put(24.6000,-4.7000){\makebox(0,0)[lb]{$M^3$}}%
\put(19.8000,-7.3000){\makebox(0,0)[lb]{$M^2$}}%
\put(3.2000,-17.3000){\makebox(0,0)[lt]{{\small The mean curvature flows starting from leaves $M^1,M^2$ and $M^3$ of}}}%
%
\special{pn 20}%
\special{sh 1}%
\special{ar 2180 1300 10 10 0  6.28318530717959E+0000}%
\special{sh 1}%
\special{ar 2180 1300 10 10 0  6.28318530717959E+0000}%
%
\special{pn 13}%
\special{pa 1990 1470}%
\special{pa 2130 1350}%
\special{fp}%
\special{sh 1}%
\special{pa 2130 1350}%
\special{pa 2066 1378}%
\special{pa 2090 1386}%
\special{pa 2092 1410}%
\special{pa 2130 1350}%
\special{fp}%
%
\special{pn 13}%
\special{pa 1810 1160}%
\special{pa 2110 1270}%
\special{fp}%
\special{sh 1}%
\special{pa 2110 1270}%
\special{pa 2054 1228}%
\special{pa 2060 1252}%
\special{pa 2042 1266}%
\special{pa 2110 1270}%
\special{fp}%
%
\special{pn 13}%
\special{pa 2590 1230}%
\special{pa 2250 1290}%
\special{fp}%
\special{sh 1}%
\special{pa 2250 1290}%
\special{pa 2320 1298}%
\special{pa 2304 1282}%
\special{pa 2312 1260}%
\special{pa 2250 1290}%
\special{fp}%
\put(22.1000,-12.4000){\makebox(0,0)[lb]{$e$}}%
\put(3.2000,-19.1000){\makebox(0,0)[lt]{{\small ${\mathfrak F}_{\{0\},\bar{\it l}_1,\cdots,\bar{\it l}_k}$ converge to the only minimal leaf $M^0$ of ${\mathfrak F}_{\{0\},\bar{\it l}_1,\cdots,\bar{\it l}_k}$.}}}%
%
\special{pn 8}%
\special{ar 2260 1300 80 840  3.1566600 4.3205356}%
\put(21.6000,-5.0000){\makebox(0,0)[lb]{$M^0$}}%
\end{picture}%
\hspace{2truecm}}

\vspace{0.4truecm}

\centerline{{\bf Figure 3.}}

\vspace{0.4truecm}

The following question arises naturally.  

\vspace{0.4truecm}

\noindent
{\it Question.} {\sl Let $\mathfrak F$ be an isoparametric foliation consisting of non-compact regular leaves 
on a non-positively curved Riemannian manifold.  Assume that the leaves of $\mathfrak F$ are 
cohomogeneity compact (i.e., each leaf $L$ is invariant under some subgroup action $H_L$ of the isometry group 
of the ambient space and the quotient space $L/H_L$ is compact).  In what case, does the result similar to 
Theorem A or B hold for $\mathfrak F$?}

\section{Mean curvature flow.}
In this section, we shall recall the notion of the mean curvature flow.  
Let $f_t$'s ($t\in[0,T)$) be a one-parameter $C^{\infty}$-family of immersions 
of a manifold $M$ into a Riemannian manifold $\widetilde M$, where $T$ is a positive constant or $T=\infty$.  
Define a map $F:M\times[0,T)\to\widetilde M$ by $F(x,t)=f_t(x)$ ($(x,t)\in M\times[0,T)$).  
Denote by $\pi$ the natural projection of $M\times[0,T)$ onto $M$.  
For a vector bundle $E$ over $M$, denote by $\pi^{\ast}E$ the induced bundle of $E$ by $\pi$.  
Also, denote by $H_t$ and $g_t$ the mean curvature vector field and the induced metric of $f_t$, respectively.  
Define a section $g$ of $\pi^{\ast}(T^{(0,2)}M)$ by 
$g_{(x,t)}:=(g_t)_x$ ($(x,t)\in M\times [0,T)$) and sections $H$ of 
$F^{\ast}T\widetilde M$ by $H_{(x,t)}:=(H_t)_x$ ($(x,t)\in M\times [0,T)$), where 
$T^{(0,2)}M$ is the tensor bundle of degree $(0,2)$ of $M$ and $T\widetilde M$ is the tangent bundle of 
$\widetilde M$.  The family $f_t$'s ($0\leq t<T$) is called a {\it mean curvature flow} if it satisfies 
$$F_{\ast}\left(\frac{\partial}{\partial t}\right)=H.\leqno{(1.1)}$$
In particular, if $f_t$'s are embeddings, then we call $M_t:=f_t(M)$'s 
$(0\in[0,T))$ rather than $f_t$'s ($0\leq t<T$) a mean curvature flow.  
See [H1,2] and [B] and so on about the study of the mean curvature flow 
(treated as the evolution of an immersion).  

\section{The non-curvature-adaptedness of the leaves.}
In [K1], we proved the following statement:

\vspace{0.3truecm}

\noindent
($\ast$) $\,\,$ {\sl 
If the root system $\triangle$ of $G/K$ is non-reduced and 
$2\lambda_{i_0}\in\triangle_+$ for some $i_0\in\{1,\cdots,k\}$, then all leaves of 
${\mathfrak F}_{\mathfrak b,\overline{\it l}_1,\cdots,\overline{\it l}_k}$ are not curvature-adapted.}

\vspace{0.3truecm}

\noindent
(see the statement (ii) of Proposition 3.5 in [K1]).  
However, the second-half part of the proof was false.  
Hence we shall correct the proof of this statement by recalculating the normal Jacobi operators of 
the leaves (see Proposition 3.5).  We shall use the notations in Introduction.  
According to the fact (iv) stated in Introduction, we have 
$$L_{b\cdot\gamma_{\xi^1}(t_1)\cdot\,\cdots\,\cdot\gamma_{\xi^k}(t_k)}
(S_{\mathfrak b,{\it l}_1,\cdots,{\it l}_k}\cdot e)
=S_{\mathfrak b,\overline{\it l}_1,\cdots,\overline{\it l}_k}\cdot
(b\cdot\gamma_{\xi^1}(t_1)\cdot\,\cdots\,\cdot\gamma_{\xi^k}(t_k)).$$
Hence we suffice to show that the leaves $S_{\mathfrak b,{\it l}_1,\cdots,{\it l}_k}\cdot e$'s are not 
curvature-adapted.  
As stated in Example 2, we set $\displaystyle{\xi^i_{t_i}:=\frac{1}{\cosh(\vert\vert\lambda_i\vert\vert t_i)}
\xi^i-\frac{1}{\vert\vert\lambda_i\vert\vert}\tanh(\vert\vert\lambda_i\vert\vert t_i)H_{\lambda_i}}$.  
For the shape operator of $S_{\mathfrak b,{\it l}_1,\cdots,{\it l}_k}\cdot e$, 
we showed the following facts (see Lemma 3.2 of [K1]).  

\vspace{0.5truecm}

\noindent
{\bf Lemma 3.1([K1]).} {\sl 
Let $A$ be the shape tensor of 
$S_{\mathfrak b,{\it l}_1,\cdots,{\it l}_k}\cdot e\,(\subset AN)$.  
Then, for $A_{\xi_0}$ ($\xi_0\in\mathfrak b$) and $A_{\xi^i_{t_i}}$ ($i=1,\cdots,k$), 
the following statements ${\rm(i)}\sim{\rm(vii)}$ hold:

{\rm(i)} For $X\in\mathfrak a\ominus(\mathfrak b+\sum\limits_{i=1}^k{\rm R}
H_{\lambda_i})$, we have $A_{\xi_0}X=A_{\xi^i_{t_i}}X=0$ 
($i=1,\cdots,k$).

{\rm (ii)} For $X\in{\rm Ker}({\rm ad}(\xi^i)\vert_{\mathfrak g_{\lambda_i}})
\ominus{\bf R}\xi^i$, 
we have $A_{\xi_0}X=0$ and $A_{\xi^i_{t_i}}X=-\vert\vert\lambda_i\vert\vert
\tanh(\vert\vert\lambda_i\vert\vert t_i)X$.  

{\rm (iii)} Assume that $2\lambda_i\in\triangle_+$.  For 
$X\in\mathfrak g_{2\lambda_i}$, we have 
$A_{\xi_0}([\theta\xi^i,X])=0$ and 
$$\begin{array}{l}
\displaystyle{A_{\xi^i_{t_i}}X=-2\vert\vert\lambda_i\vert\vert\tanh(
\vert\vert\lambda_i\vert\vert t_i)X-\frac{1}{2\cosh(\vert\vert\lambda_i\vert\vert t_i)}
[\theta\xi^i,X],}\\
\displaystyle{A_{\xi^i_{t_i}}([\theta\xi^i,X])
=-\frac{\vert\vert\lambda_i\vert\vert^2}{\cosh(\vert\vert\lambda_i\vert\vert t_i)}X
-\vert\vert\lambda_i\vert\vert\tanh(\vert\vert\lambda_i\vert\vert t_i)[\theta\xi^i,X],}
\end{array}$$
where $\theta$ is the Cartan involution of $\mathfrak g$ with 
${\rm Fix}\,\theta=\mathfrak k$.  

{\rm(iv)} For $X\in({\bf R}\xi^i+{\bf R}H_{\lambda_i})\ominus{\it l}_i$, 
we have $A_{\xi_0}X=0$ and $A_{\xi^i_{t_i}}X=-\vert\vert\lambda_i\vert\vert
\tanh(\vert\vert\lambda_i\vert\vert t_i)X$.

{\rm (v)} For $X\in(\mathfrak g_{\lambda_j}\ominus{\bf R}\xi^j)
+(({\bf R}\xi^j+{\bf R}H_{\lambda_j})\ominus{\it l}_j)
+\mathfrak g_{2\lambda_j}$ ($j\not=i$), we have 
$A_{\xi_0}X=A_{\xi^i_{t_i}}X=0$.

{\rm(vi)} For $X\in\mathfrak g_{\mu}$ ($\mu\in\triangle_+\setminus
\{\lambda_1,\cdots,\lambda_k\}$), we have $A_{\xi_0}X=\mu(\xi_0)X$.

{\rm (vii)} Let 
$\displaystyle{k_i:=\exp\left(\frac{\pi}{\sqrt2\vert\vert\lambda_i\vert\vert}
(\xi^i+\theta\xi^i)\right)}$, where $\exp$ is the exponential map of $G$.  
Then ${\rm Ad}(k_i)\circ A_{\xi^i_{t_i}}=-A_{\xi^i_{t_i}}\circ{\rm Ad}(k_i)$ 
holds over $\mathfrak n\ominus\sum\limits_{i=1}^k(\mathfrak g_{\lambda_i}
+\mathfrak g_{2\lambda_i})$, where ${\rm Ad}$ is the adjoint representation 
of $G$.}

\vspace{0.5truecm}

\noindent
{\it Remark 3.1.} If $\lambda_i\in\triangle_+$, then we have $\vert\vert\lambda_i\vert\vert=\sqrt 2$ 
from how to choose the metric of $G/K$ (see Introduction).  

\vspace{0.5truecm}

According to (5.3) in Page 310 of [M], we have the following fact.  

\vspace{0.5truecm}

\noindent
{\bf Lemma 3.2([M]).} {\sl Let $X$ and $Y$ be left-invariant vector fields on $AN$ and $\nabla$ be 
the Levi-Civita connection of the left-invariant metric 
$\langle\,\,,\,\,\rangle$ of $AN$.  Then we have 
$$\begin{array}{l}
\displaystyle{\nabla_XY=\frac{1}{2}\left(\,\,[X,Y]-{\rm ad}(X)^{\ast}(Y)-{\rm ad}(Y)^{\ast}(X)\,\,\right),}
\end{array}\leqno{(3.2)}$$
where ${\rm ad}(X)^{\ast}$ (resp. ${\rm ad}(Y)^{\ast}$) is the adjoint operator of ${\rm ad}(X)$ 
(resp. ${\rm ad}(Y)$) with respect to $\langle\,\,,\,\,\rangle_e$ and 
$(\bullet)_{\mathfrak a+\mathfrak n}$ is the the $(\mathfrak a+\mathfrak n)$-component of $(\bullet)$.}

\vspace{0.5truecm}

Let ${\rm pr}^1_{\mathfrak a+\mathfrak n}$ (resp. 
${\rm pr}^2_{\mathfrak a+\mathfrak n}$) be the projection of $\mathfrak g$ 
onto $\mathfrak a+\mathfrak n$ with respect to the decomposition 
$\mathfrak g=\mathfrak k+(\mathfrak a+\mathfrak n)$ (resp. 
$\mathfrak g=(\mathfrak k_0+\sum\limits_{\lambda\in\triangle_+}
\mathfrak p_{\lambda})+(\mathfrak a+\mathfrak n)$).  
We ([K1]) showed the following facts (see the proof of Lemma 3.2 in [K1]).  

\vspace{0.5truecm}

\noindent
{\bf Lemma 3.3([K1]).} {\sl {\rm (i)} For any $H\in\mathfrak a$, we have 
$${\rm ad}(H)^{\ast}={\rm ad}(H).\leqno{(3.3)}$$

{\rm (ii)} For any $X\in\mathfrak g_{\lambda}$, we have 
$$\begin{array}{l}
\displaystyle{{\rm ad}(X)^{\ast}=-{\rm pr}_{\mathfrak a+\mathfrak n}\circ{\rm ad}(\theta X)}\\
\hspace{1.35truecm}\displaystyle{=\left\{
\begin{array}{cl}
\displaystyle{0} & \displaystyle{{\rm on}\,\,\mathfrak a}\\
\displaystyle{
\begin{array}{l}
\displaystyle{-\langle X,\cdot\rangle_e\otimes H_{\lambda}-
{\rm pr}_{\mathfrak n}\circ{\rm pr}^1_{\mathfrak a+\mathfrak n}\circ
{\rm ad}(X_{\mathfrak k})}\\
\displaystyle{+{\rm pr}_{\mathfrak n}\circ{\rm pr}^2_{\mathfrak a+\mathfrak n}
\circ{\rm ad}(X_{\mathfrak p})}
\end{array}
} & \displaystyle{{\rm on}\,\,\mathfrak n},
\end{array}
\right.}
\end{array}
\leqno{(3.4)}$$
where $(\bullet)_{\mathfrak k}$ (resp. $(\cdot)_{\mathfrak p}$) denotes the $\mathfrak k$-component 
(resp. $\mathfrak p$-component) of $(\bullet)$.}

\vspace{0.5truecm}

According to $(3.4)$, we have 
$${\rm ad}(X)^{\ast}(Y)=\left\{
\begin{array}{ll}
0 & (\lambda-\mu\in\triangle_+)\\
-\langle X,Y\rangle H_{\lambda} & (\lambda=\mu)\\
-[\theta X,Y] & (\mu-\lambda\in\triangle_+)\\
0 & (\lambda-\mu\notin\triangle\cup\{0\})
\end{array}\right.
\leqno{(3.5)}$$
for any $X\in\mathfrak g_{\lambda}$ ($\lambda\in\triangle_+$) and any 
$Y\in\mathfrak g_{\mu}$ ($\mu\in\triangle_+$).  
For each $X\in\mathfrak a+\mathfrak n$, we denote by $\widetilde X$ the left-invariant vector field on $AN$ with 
$(\widetilde X)_e=X$.  
By using Lemma 3.2, $(3.3),\,(3.4)$ and $(3.5)$, we can derive the facts directly.  

\vspace{0.5truecm}

\noindent
{\bf Lemma 3.4.} {\sl For any unit vector $X_{\lambda},\,Y_{\lambda}$ of $\mathfrak g_{\lambda}$ 
($\lambda\in\triangle_+$) and $H_{\lambda}$ ($\lambda\in\triangle_+$), we have 
$$
\nabla_{\widetilde H_{\lambda}}\widetilde H_{\mu}=\nabla_{\widetilde H_{\lambda}}\widetilde X_{\mu}
=0,\,\,\,\,
\nabla_{\widetilde X_{\lambda}}\widetilde H_{\mu}=-\lambda(H_{\mu})\widetilde X_{\lambda}
\quad(\lambda,\,\,\mu\in\triangle_+)$$
and 
$$\nabla_{\widetilde X_{\lambda}}\widetilde Y_{\mu}=\left\{
\begin{array}{ll}
\displaystyle{\frac{1}{2}\left([\widetilde X_{\lambda},\widetilde Y_{\mu}]
+\widetilde{\theta[Y_{\mu},\theta} X_{\lambda}]\right)} & 
(\lambda-\mu\in\triangle_+)\\
\displaystyle{\frac{1}{2}[\widetilde X_{\lambda},\widetilde Y_{\mu}]
+\langle\widetilde X_{\lambda},\widetilde Y_{\mu}\rangle\widetilde H_{\lambda}} & 
(\lambda=\mu)\\
\displaystyle{\frac{1}{2}\left([\widetilde X_{\lambda},\widetilde Y_{\mu}]
+\widetilde{\theta[X_{\lambda},\theta} Y_{\mu}]\right)} & 
(\mu-\lambda\in\triangle_+)\\
\displaystyle{\frac{1}{2}[\widetilde X_{\lambda},\widetilde Y_{\mu}]} & (\lambda-\mu\notin\triangle\cup\{0\})
\end{array}
\right.$$
}

\vspace{0.5truecm}

From Lemma 3.4 and $(3.5)$, we can derive the following facts for the normal Jacobi operators 
by somewhat long calculations.  

\vspace{0.5truecm}

\noindent
{\bf Proposition 3.5.} {\sl Let $R$ be the curvature tensor of $AN(=G/K)$.  
Then, for $R(\xi_0)$ ($\xi^0\in\mathfrak b$) and $R(\xi^i_{t_i})$ ($i=1,\cdots,k$), 
the following statements ${\rm(i)}\sim{\rm(vi)}$ hold:

{\rm(i)} For $X\in\mathfrak a\ominus(\mathfrak b+\sum\limits_{i=1}^k{\rm R}
H_{\lambda_i})$, we have $R(\xi_0)(X)=R(\xi^i_{t_i})(X)=0$ 
($i=1,\cdots,k$).

{\rm (ii)} For $X\in{\rm Ker}({\rm ad}(\xi^i)\vert_{\mathfrak g_{\lambda_i}})
\ominus{\bf R}\xi^i$, 
we have $R(\xi_0)(X)=0$ and $R(\xi^i_{t_i})(X)
=\frac{\vert\vert\lambda_i\vert\vert^2}{2}(1-3\tanh^2(\vert\vert\lambda_i\vert\vert t_i))X$.  

{\rm (iii)} Assume that $2\lambda_i\in\triangle_+$ (hence $\vert\vert\lambda_i\vert\vert=\sqrt 2$).  
For $X\in\mathfrak g_{2\lambda_i}$, we have 
$R(\xi_0)(X)=R(\xi_0)([\theta\xi^i,X])=0$ and 
$$\begin{array}{l}
\displaystyle{R(\xi^i_{t_i})(X)=-\vert\vert\lambda_i\vert\vert^2(1+3\tanh^2(\vert\vert\lambda_i\vert\vert t_i))X
-\frac{3\vert\vert\lambda_i\vert\vert\tanh(\vert\vert\lambda_i\vert\vert t_i)}
{2\cosh(\vert\vert\lambda_i\vert\vert t_i)}[\theta\xi^i,X]
}\\
\displaystyle{R(\xi^i_{t_i})([\theta\xi^i,X])
=-\frac{6\vert\vert\lambda_i\vert\vert\tanh(\vert\vert\lambda_i\vert\vert t_i)}
{\cosh(\vert\vert\lambda_i\vert\vert t_i)}X
+\frac{\sqrt 2\vert\vert\lambda_i\vert\vert}{4}(1-3\tanh^2(\vert\vert\lambda_i\vert\vert t_i))
[\theta\xi^i,X].}
\end{array}$$

{\rm(iv)} For $X\in({\bf R}\xi^i+{\bf R}H_{\lambda_i})\ominus{\it l}_i$, 
we have $R(\xi_0)(X)=0$ and $R(\xi^i_{t_i})(X)=-\vert\vert\lambda_i\vert\vert^2X$.  

{\rm (v)} For $X\in(\mathfrak g_{\lambda_j}\ominus{\bf R}\xi^j)
+(({\bf R}\xi^j+{\bf R}H_{\lambda_j})\ominus{\it l}_j)
+\mathfrak g_{2\lambda_j}$ ($j\not=i$), we have 
$R(\xi_0)(X)=R(\xi^i_{t_i})(X)=0$.  

{\rm(vi)} For $X\in\mathfrak g_{\mu}$ ($\mu\in\triangle_+\setminus
\{\lambda_1,\cdots,\lambda_k\}$), we have $R(\xi_0)(X)=-\mu(\xi_0)^2X$.}

\vspace{0.5truecm}

From Lemma 3.1 and Proposition 3.5, we can derive the following facts directly.  

\vspace{0.5truecm}

\noindent
{\bf Proposition 3.6.} {\sl 
For $[A_{\xi_0},R(\xi_0)]$ ($\xi_0\in\mathfrak b$) and $[A_{\xi^i_{t_i}},R(\xi^i_{t_i})]$ 
($i=1,\cdots,k$), the following statements ${\rm(i)}\sim{\rm(vi)}$ hold:

{\rm(i)} For $X\in\mathfrak a\ominus(\mathfrak b+\sum\limits_{i=1}^k{\rm R}H_{\lambda_i})$, we have 
$[A,R(\xi_0)](X)=[A_{\xi^i_{t_i}},R(\xi^i_{t_i})](X)=0$ ($i=1,\cdots,k$).  

{\rm (ii)} For $X\in{\rm Ker}({\rm ad}(\xi^i)\vert_{\mathfrak g_{\lambda_i}})\ominus{\bf R}\xi^i$, 
we have $[A_{\xi_0},R(\xi_0)](X)=[A_{\xi^i_{t_i}},R(\xi^i_{t_i})](X)=0$.  

{\rm (iii)} Assume that $2\lambda_i\in\triangle_+$ (hence $\vert\vert\lambda_i\vert\vert=\sqrt 2$).  
For $X\in\mathfrak g_{2\lambda_i}$, we have 
$[A_{\xi_0},R(\xi_0)](X)=[A_{\xi_0},R(\xi_0)]([\theta\xi^i,X])=0$ and 
$$\begin{array}{l}
\displaystyle{[A_{\xi^i_{t_i}},R(\xi^i_{t_i})](X)
=-\frac{3}{2\cosh^3(\sqrt 2t_i)}[\theta\xi^i,X]}\\
\displaystyle{[A_{\xi^i_{t_i}},R(\xi^i_{t_i})]([\theta\xi^i,X])
=-\frac{6}{\cosh^3(\sqrt 2t_i)}X.}
\end{array}$$

{\rm(iv)} For $X\in({\bf R}\xi^i+{\bf R}H_{\lambda_i})\ominus{\it l}_i$, 
we have $[A_{\xi_0},R(\xi_0)](X)=[A_{\xi^i_{t_i}},R(\xi^i_{t_i})](X)=0$.  

{\rm (v)} For $X\in(\mathfrak g_{\lambda_j}\ominus{\bf R}\xi^j)
+(({\bf R}\xi^j+{\bf R}H_{\lambda_j})\ominus{\it l}_j)
+\mathfrak g_{2\lambda_j}$ ($j\not=i$), we have 
$[A_{\xi_0},R(\xi_0)](X)=[A_{\xi^i_{t_i}},R(\xi^i_{t_i})](X)=0$.  

{\rm(vi)} For $X\in\mathfrak g_{\mu}$ ($\mu\in\triangle_+\setminus
\{\lambda_1,\cdots,\lambda_k\}$), we have 
$[A_{\xi_0},R(\xi_0)](X)=[A_{\xi^i_{t_i}},R(\xi^i_{t_i})](X)=0$.}

\vspace{0.5truecm}

From (iv) of Proposition 3.6, we can derive the statement $(\ast)$.  

\vspace{0.5truecm}

Also, we ([K1]) showed the following fact in terms of Lemma 3.1.  

\vspace{0.5truecm}

\noindent
{\bf Proposition 3.7([K1]).} {\sl If $\mathfrak b=\{0\}$, then 
${\mathfrak F}_{\mathfrak b,{\it l}_{\xi^1,t_1},\cdots,{\it l}_{\xi^k,t_k}}$ admits the only minimal leaf.}

\section{Proof of Theorem A} 
In this section, we shall prove Theorem A.  
We use the notations in Sections 1 and 3.  
Note that ${\rm Exp}\vert_{\mathfrak a}=\exp\vert_{\mathfrak a}$ and 
${\rm Exp}\vert_{\mathfrak n}\not=\exp\vert_{\mathfrak n}$.  
Set 
$\Sigma:={\rm Exp}(T^{\perp}_eS_{\mathfrak b,\overline{\it l}_1,\cdots,\overline{\it l}_k}\cdot e)
(={\rm Exp}(\mathfrak b+{\mathbb R}\{\xi^1,\cdots,\xi^k\}))$, 
which is the flat section of the $S_{\mathfrak b,\overline{\it l}_1,\cdots,\overline{\it l}_k}$-action through 
$e$.  Each leaf of ${\mathfrak F}_{\mathfrak b,\overline{\it l}_1,\cdots,\overline{\it l}_k}$ meets $\Sigma$ 
at the only one point.  That is, $\Sigma$ is regarded as the leaf space of this foliation.  
For $\xi_0\in\mathfrak b$ and $t_i\in\mathbb R\,\,(i=1,\cdots,k)$, we set 
$x_{\xi_0,t_1,\cdots,t_k}:={\rm Exp}\xi_0\cdot\gamma_{\xi^1(t_1)}\cdot\,\cdots\,\cdot\gamma_{\xi^k(t_k)}$.  
Also, denote by $\frac{D}{ds}(\bullet)$ the covariant derivative of vector fields $(\bullet)$ along curves in 
$AN$ (with respect to the left-invariant metric).  
The following fact is well-known about the geodesics in rank one symmetric spaces of non-compact type 
but we shall give the proof.  

\vspace{0.5truecm}

\noindent
{\bf Lemma 4.1.} {\sl The velocity vector $\gamma_{\xi^i}'(s)$ ($i=1,\cdots,k$) 
is described as 
$$\gamma_{\xi^i}'(s)=\frac{1}{\cosh(\vert\vert\lambda_i\vert\vert s)}
(\widetilde{\xi^i})_{{\gamma}_{\xi^i}(s)}-\frac{\tanh(\vert\vert\lambda_i\vert\vert s)}
{\vert\vert\lambda_i\vert\vert}(\widetilde{H_{\lambda_i}})_{\gamma_{\xi^i}(s)}\leqno{(4.1)}$$
and $\gamma_{\xi_0}'(s)$ is described as 
$$\gamma_{\xi_0}'(s)=(\widetilde{\xi}_0)_{\gamma_{\xi_0}(s)}\leqno{(4.2)}.$$}

\vspace{0.5truecm}

\noindent
{\it Proof.} Set $Y(s):=\frac{1}{\cosh(\vert\vert\lambda_i\vert\vert s)}
(\widetilde{\xi^i})_{{\gamma}_{\xi^i}(s)}-\frac{\tanh(\vert\vert\lambda_i\vert\vert s)}
{\vert\vert\lambda_i\vert\vert}(\widetilde{H_{\lambda_i}})_{\gamma_{\xi^i}(s)}$.  
It is clear that $Y(0)=\xi^i$.  By using Lemma 3.4, we can show $\frac{D}{ds}Y=0$.  
Hence we obtain $Y(s)=\gamma_{\xi^i}'(s)$.  
Also, it is clear that $(\widetilde{\xi}_0)_{\gamma_{\xi_0}(0)}=\xi_0$.  By using Lemma 3.4, we can show 
$\frac{D}{ds}(\widetilde{\xi}_0)_{\gamma_{\xi_0}(s)}=0$.  
Hence we obtain $(\widetilde{\xi}_0)_{\gamma_{\xi_0}(s)}=\gamma_{\xi_0}'(s)$.  \hspace{9.5truecm}q.e.d.

\vspace{0.5truecm}

Next we shall show the following fact.  

\vspace{0.5truecm}

\noindent
{\bf Lemma 4.2.} {\sl The point $x_{\xi_0,t_1,\cdots,t_k}$ belongs to $\Sigma$.}

\vspace{0.5truecm}

\noindent
{\it Proof.} It is clear that ${\rm Exp}(\xi_0)$ belongs to $\Sigma$.  
First we shall show that ${\rm Exp}(\xi_0)\cdot\gamma_{\xi^1(t_1)}$ belongs to $\Sigma$.  
Let $\gamma_{\xi_0}$ be the geodesic in $AN$ with $\gamma_{\xi_0}'(0)=\xi_0$.  
Since $\gamma_{\xi^1}$ is a geodesic in $AN$ and $L_{{\rm Exp}(\xi_0)}$ is an isometry of $AN$, 
$L_{{\rm Exp}(\xi_0)}\circ\gamma_{\xi^1}$ is a geodesic in $AN$.  
Hence we suffice to show that $(L_{{\rm Exp}(\xi_0)}\circ\gamma_{\xi^1})'(0)
=(\widetilde{\xi}^1)_{{\rm Exp}(\xi_0)}$ is tangent to $\Sigma$.  
Denote by $\widehat{\xi}^1$ the parallel vector field along $\gamma_{\xi_0}$.  
Take orthonormal bases $\{e^{\lambda}_1,\cdots,e^{\lambda}_{m_{\lambda}}\}$ of $\mathfrak g_{\lambda}$ 
($\lambda\in\triangle_+$).  Also, take an orthonormal base $\{e^0_1,\cdots,e^0_r\}$ of $\mathfrak a$.  
We describe $\widehat{\xi}^1$ as 
$$
\widehat{\xi}^1(s)=\sum_{i=1}^ra^0_i(s)(\widetilde{e^0_i})_{\gamma_{\xi_0}(s)}
+\sum_{\lambda\in\triangle_+}\sum_{i=1}^{m_{\lambda}}
a^{\lambda}_i(s)(\widetilde{e^{\lambda}_i})_{\gamma_{\xi_0}(s)}\quad(s\in{\mathbb R}),$$
where $a^0_i$ and $a^{\lambda}_i$ are functions over ${\mathbb R}$.  
Fix $s_0\in\mathbb R$.  By using Lemma 3.4, we can show 
$$\begin{array}{l}
\displaystyle{\left.\frac{D}{ds}\right\vert_{s=s_0}\widehat{\xi}^1
=\sum_{i=1}^r\left((a^0_i)'(s_0)(\widetilde{e^0_i})_{\gamma_{\xi_0}(s_0)}
+(a^0_i)(s_0)\left.\frac{D}{ds}\right\vert_{s=s_0}((\widetilde{e^0_i})_{\gamma_{\xi_0}(s)})\right)}\\
\hspace{0.5truecm}\displaystyle{
+\sum_{\lambda\in\triangle_+}\sum_{i=1}^{m_{\lambda}}
\left((a^{\lambda}_i)'(s_0)(\widetilde{e^{\lambda}_i})_{\gamma_{\xi_0}(s_0)}
+a^{\lambda}_i(s_0)\left.\frac{D}{ds}\right\vert_{s=s_0}
((\widetilde{e^{\lambda}_i})_{\gamma_{\xi_0}(s)})\right)}\\
\displaystyle{=\sum_{i=1}^r\left((a^0_i)'(s_0)(\widetilde{e^0_i})_{\gamma_{\xi_0}(s_0)}
+(a^0_i)(s_0)\nabla_{\gamma_{\xi_0}'(s_0)}((\widetilde{e^0_i})_{\gamma_{\xi_0}(s_0)})\right)}\\
\hspace{0.5truecm}\displaystyle{
+\sum_{\lambda\in\triangle_+}\sum_{i=1}^{m_{\lambda}}
\left((a^{\lambda}_i)'(s_0)(\widetilde{e^{\lambda}_i})_{\gamma_{\xi_0}(s_0)}
+a^{\lambda}_i(s_0)\nabla_{\gamma_{\xi_0}'(s_0)}((\widetilde{e^{\lambda}_i})_{\gamma_{\xi_0}(s_0)})\right)}\\
\displaystyle{=\sum_{i=1}^r\left((a^0_i)'(s_0)(\widetilde{e^0_i})_{\gamma_{\xi_0}(s_0)}
+(a^0_i)(s_0)(\nabla_{\widetilde{\xi_0}}\widetilde{e^0_i})_{\gamma_{\xi_0}(s_0)}\right)}\\
\hspace{0.5truecm}\displaystyle{
+\sum_{\lambda\in\triangle_+}\sum_{i=1}^{m_{\lambda}}
\left((a^{\lambda}_i)'(s_0)(\widetilde{e^{\lambda}_i})_{\gamma_{\xi_0}(s_0)}
+a^{\lambda}_i(s_0)(\nabla_{\widetilde{\xi_0}}\widetilde{e^{\lambda}_i})_{\gamma_{\xi_0}(s_0)})\right)}\\
\displaystyle{=\sum_{i=1}^r(a^0_i)'(s_0)(\widetilde{e^0_i})_{\gamma_{\xi_0}(s_0)}
+\sum_{\lambda\in\triangle_+}\sum_{i=1}^{m_{\lambda}}
(a^{\lambda}_i)'(s_0)(\widetilde{e^{\lambda}_i})_{\gamma_{\xi_0}(s_0)}=0,}
\end{array}$$
that is, $(a^0_i)'(s_0)=(a^{\lambda}_i)'(s_0)=0$, 
where we use $\gamma_{\xi_0}'(s_0)=\widetilde{\xi_0}_{\gamma_{\xi_0}(s_0)}$.  
From the arbitrariness of $s_0$, we see that $a^0_i$ and $a^{\lambda}_i$ are constant.  
Hence we obtain $\widehat{\xi}^1(s)=(\widetilde{\xi}^1)_{\gamma_{\xi_0}(s)}$.  
On the other hand, since $\xi^1$ is tangent to $\Sigma$ and $\Sigma$ is totally geodesic, $\widehat{\xi}^1(1)$ 
also is tangent to $\Sigma$.  
Hence we see that $(\widetilde{\xi}^1)_{{\rm Exp}(\xi_0)}$ is tangent to $\Sigma$.  
Therefore ${\rm Exp}(\xi_0)\cdot\gamma_{\xi^1(t_1)}$ belongs to $\Sigma$.  

Next we shall show that ${\rm Exp}(\xi_0)\cdot\gamma_{\xi^1(t_1)}\cdot\gamma_{\xi^2(t_2)}$ belongs to $\Sigma$.  
Since $\gamma_{\xi^2}$ is a geodesic in $AN$ and $L_{{\rm Exp}(\xi_0)\cdot\gamma_{\xi^1}(t_1)}$ is an isometry of 
$AN$, $L_{{\rm Exp}(\xi_0)\cdot\gamma_{\xi^1}(t_1)}\circ\gamma_{\xi^2}$ is a geodesic in $AN$.  
Hence we suffice to show that $(L_{{\rm Exp}(\xi_0)\cdot\gamma_{\xi^1}(t_1)}\circ\gamma_{\xi^2})'(0)
=(\widetilde{\xi}^2)_{{\rm Exp}(\xi_0)\cdot\gamma_{\xi^1}(t_1)}$ is tangent to $\Sigma$.  
Denote by $\widehat{\xi}^2$ the parallel vector field along 
$\overline{\gamma}_{\xi_1}:=L_{{\rm Exp}(\xi_0)}\circ\gamma_{\xi^1}$ with 
$\widehat{\xi}^2(0)=(\widetilde{\xi}^2)_{{\rm Exp}(\xi_0)}$.  
We describe $\widehat{\xi}^2$ as 
$$
\widehat{\xi}^2(s)=\sum_{i=1}^rb^0_i(s)(\widetilde{e^0_i})_{\overline{\gamma}_{\xi^1}(s)}
+\sum_{\lambda\in\triangle_+}\sum_{i=1}^{m_{\lambda}}
b^{\lambda}_i(s)(\widetilde{e^{\lambda}_i})_{\overline{\gamma}_{\xi^1}(s)}\quad(s\in{\mathbb R}),$$
where $b^0_i$ and $b^{\lambda}_i$ are functions over ${\mathbb R}$.  
Fix $s_0\in\mathbb R$.  By using Lemma 3.4, we can show 
$$\begin{array}{l}
\displaystyle{\left.\frac{D}{ds}\right\vert_{s=s_0}\widehat{\xi}^2
=\sum_{i=1}^r\left((b^0_i)'(s_0)(\widetilde{e^0_i})_{\overline{\gamma}_{\xi^1}(s_0)}
+(b^0_i)(s_0)\left.\frac{D}{ds}\right\vert_{s=s_0}((\widetilde{e^0_i})_{\overline{\gamma}_{\xi^1}(s)})\right)}\\
\hspace{0.5truecm}\displaystyle{
+\sum_{\lambda\in\triangle_+}\sum_{i=1}^{m_{\lambda}}
\left((b^{\lambda}_i)'(s_0)(\widetilde{e^{\lambda}_i})_{\overline{\gamma}_{\xi^1}(s_0)}
+b^{\lambda}_i(s_0)\left.\frac{D}{ds}\right\vert_{s=s_0}
((\widetilde{e^{\lambda}_i})_{\overline{\gamma}_{\xi^1}(s)})\right)}\\
\displaystyle{=\sum_{i=1}^r\left((b^0_i)'(s_0)(\widetilde{e^0_i})_{\overline{\gamma}_{\xi^1}(s_0)}
+(b^0_i)(s_0)\nabla_{\overline{\gamma}_{\xi^1}'(s_0)}((\widetilde{e^0_i})_{\overline{\gamma}_{\xi^1}(s)})\right)
}\\
\hspace{0.5truecm}\displaystyle{
+\sum_{\lambda\in\triangle_+}\sum_{i=1}^{m_{\lambda}}
\left((b^{\lambda}_i)'(s_0)(\widetilde{e^{\lambda}_i})_{\overline{\gamma}_{\xi^1}(s_0)}
+b^{\lambda}_i(s_0)\nabla_{\overline{\gamma}'_{\xi^1}(s_0)}
((\widetilde{e^{\lambda}_i})_{\overline{\gamma}_{\xi^1}(s)})\right)=0.}
\end{array}\leqno{(4.3)}$$
Since $\gamma'_{\xi^1}(s_0)=\frac{1}{\cosh(\vert\vert\lambda_1\vert\vert s_0)}
(\widetilde{\xi^1})_{{\gamma}_{\xi^1}(s_0)}-\frac{\tanh(\vert\vert\lambda_1\vert\vert s_0)}
{\vert\vert\lambda_1\vert\vert}(\widetilde{H_{\lambda_1}})_{\gamma_{\xi^1}(s_0)}$ by Lemma 4.1, 
$\overline{\gamma}_{\xi^1}'(s_0)$ is described as 
$$\begin{array}{l}
\displaystyle{\overline{\gamma}_{\xi^1}'(s_0)=(L_{{\rm Exp}(\xi_0)})_{\ast}(\gamma'_{\xi^1}(s_0))}\\
\displaystyle{=\frac{1}{\cosh(\vert\vert\lambda_1\vert\vert s_0)}
(\widetilde{\xi^1})_{\overline{\gamma}_{\xi^1}(s_0)}-\frac{\tanh(\vert\vert\lambda_1\vert\vert s_0)}
{\vert\vert\lambda_1\vert\vert}(\widetilde H_{\lambda_1})_{\overline{\gamma}_{\xi^1}(s_0)}.}
\end{array}$$
Hence, by using Lemma 3.4, we have 
$$\begin{array}{l}
\displaystyle{\nabla_{\overline{\gamma}_{\xi^1}'(s_0)}((\widetilde{e^0_i})_{\overline{\gamma}_{\xi^1}}
=\frac{1}{\cosh(\vert\vert\lambda_1\vert\vert s_0)}
(\nabla_{\widetilde{\xi^1}}\widetilde{e^0_i})_{\overline{\gamma}_{\xi^1}(s_0)}
-\frac{\tanh(\vert\vert\lambda_1\vert\vert s_0)}{\vert\vert\lambda_1\vert\vert}
(\nabla_{\widetilde{H_{\lambda_1}}}\widetilde{e^0_i})_{\overline{\gamma}_{\xi^1}(s_0)}}\\
\hspace{2.8truecm}\displaystyle{=-\frac{\lambda_1(e^0_i)}{\cosh(\vert\vert\lambda_1\vert\vert s_0)}
(\widetilde{\xi^1})_{\overline{\gamma}_{\xi^1}(s_0)}}
\end{array}\leqno{(4.4)}$$
and 
$$\begin{array}{l}
\displaystyle{\nabla_{\overline{\gamma}_{\xi^1}'(s_0)}((\widetilde{e^{\lambda}_i})_{\overline{\gamma}_{\xi^1}}
=\frac{1}{\cosh(\vert\vert\lambda_1\vert\vert s_0)}
(\nabla_{\widetilde{\xi^1}}\widetilde{e^{\lambda}_i})_{\overline{\gamma}_{\xi^1}(s_0)}
-\frac{\tanh(\vert\vert\lambda_1\vert\vert s_0)}{\vert\vert\lambda_1\vert\vert}
(\nabla_{\widetilde{H_{\lambda_1}}}\widetilde{e^{\lambda}_i})_{\overline{\gamma}_{\xi^1}(s_0)}}\\
\hspace{1.8truecm}\displaystyle{=\left\{
\begin{array}{ll}
\displaystyle{\frac{1}{2\cosh(\vert\vert\lambda_1\vert\vert s_0)}
\left([\widetilde{\xi^1},\widetilde{e^{\lambda}_i}]
+\widetilde{\theta[e^{\lambda}_i,\theta} \xi^1]\right)} & 
(\lambda_1-\lambda\in\triangle_+)\\
\displaystyle{\frac{1}{2\cosh(\vert\vert\lambda_1\vert\vert s_0)}\left(
[\widetilde{\xi}^1,\widetilde{e^{\lambda}_i}]
+2\langle\widetilde{\xi^1},\widetilde{e^{\lambda}_i}\rangle\widetilde H_{\lambda_1}\right)} & 
(\lambda_1=\lambda)\\
\displaystyle{\frac{1}{2\cosh(\vert\vert\lambda_1\vert\vert s_0)}
\left([\widetilde{\xi^1},\widetilde{e^{\lambda}_i}]
+\widetilde{\theta[\xi^1,\theta} e^{\lambda}_i]\right)} & 
(\lambda-\lambda_1\in\triangle_+)\\
\displaystyle{\frac{1}{2\cosh(\vert\vert\lambda_1\vert\vert s_0)}
[\widetilde{\xi^1},\widetilde{e^{\lambda}_i}]} & (\lambda_1-\lambda\notin\triangle\cup\{0\}).
\end{array}\right.}
\end{array}\leqno{(4.5)}$$
By substituting $(4.4)$ and $(4.5)$ into $(4.3)$, we obtain 
$$\begin{array}{l}
\displaystyle{\left.\frac{D}{ds}\right\vert_{s=s_0}\widehat{\xi}^2
=\sum_{i=1}^r\left((b^0_i)'(s_0)(\widetilde{e^0_i})_{\overline{\gamma}_{\xi^1}(s_0)}
-\frac{\lambda_1(e^0_i)(b^0_i)(s_0)}{\cosh(\vert\vert\lambda_1\vert\vert s_0)}
(\widetilde{\xi^1})_{\overline{\gamma}_{\xi^1}(s_0)}\right)}\\
\hspace{2.5truecm}\displaystyle{
+\sum_{\lambda\in\triangle_+}\sum_{i=1}^{m_{\lambda}}
(b^{\lambda}_i)'(s_0)(\widetilde{e^{\lambda}_i})_{\overline{\gamma}_{\xi^1}(s_0)}}\\
\hspace{2.5truecm}\displaystyle{
+\sum_{\lambda_1-\lambda\in\triangle_+}\sum_{i=1}^{m_{\lambda}}
\frac{b^{\lambda}_i(s_0)}{2\cosh(\vert\vert\lambda_1\vert\vert s_0)}
\left([\widetilde{\xi^1},\widetilde{e^{\lambda}_i}]
+\widetilde{\theta[e^{\lambda}_i,\theta} \xi^1]\right)}\\
\hspace{2.5truecm}\displaystyle{
+\sum_{\lambda-\lambda_1\in\triangle_+}\sum_{i=1}^{m_{\lambda}}
\frac{b^{\lambda}_i(s_0)}{2\cosh(\vert\vert\lambda_1\vert\vert s_0)}
\left([\widetilde{\xi^1},\widetilde{e^{\lambda}_i}]
+\widetilde{\theta[\xi^1,\theta} e^{\lambda}_i]\right)}\\
\hspace{2.5truecm}\displaystyle{
+\sum_{\lambda-\lambda_1\notin\triangle\cup\{0\}}\sum_{i=1}^{m_{\lambda}}
\frac{b^{\lambda}_i(s_0)}{2\cosh(\vert\vert\lambda_1\vert\vert s_0)}
[\widetilde{\xi^1},\widetilde{e^{\lambda}_i}]}\\
\hspace{2.5truecm}\displaystyle{
+\sum_{i=1}^{m_{\lambda_1}}\frac{b^{\lambda_1}_i(s_0)}{2\cosh(\vert\vert\lambda_1\vert\vert s_0)}\left(
[\widetilde{\xi}^1,\widetilde{e^{\lambda_1}_i}]
+2\langle\widetilde{\xi^1},\widetilde{e^{\lambda_1}_i}\rangle\widetilde H_{\lambda_1}\right)=0.}
\end{array}\leqno{(4.6)}$$
Without loss of generality, we may assume that $e^{\lambda_2}_1=\xi^2$.  
Hence we have $b^{\lambda_2}_1(0)=1$ and $b^{\lambda}_i(0)=0$ for any $(\lambda,i)$ other than 
$(\lambda_2,1)$.  
From $(4.6)$ and these relations, we obtain 
$b^{\lambda_2}_1\equiv 1$ and $b^{\lambda}_i\equiv 0$ for any $(\lambda,i)$ other than $(\lambda_2,1)$, where 
we note that $\lambda_1-\lambda_2\notin\triangle\cup\{0\}$.  
Therefore we obtain $\widehat{\xi^2}=(\widetilde{\xi^2})_{\overline{\gamma}_{\xi^1}(s)}$.  
On the other hand, since $(\widehat{\xi^2})(0)$ is tangent to $\Sigma$ and $\Sigma$ is totally geodesic, 
$\widehat{\xi}^2(t_1)$ also is tangent to $\Sigma$.  
Hence we see that $(\widetilde{\xi^2})_{{\rm Exp}(\xi_0)\cdot\gamma_{\xi^1}(t_1)}$ is tangent to $\Sigma$.  
Therefore ${\rm Exp}(\xi_0)\cdot\gamma_{\xi^1(t_1)}\cdot\gamma_{\xi^2(t_2)}$ belongs to $\Sigma$.  
In the sequel, by repeating the same discussion, we can derive that 
$x_{\xi^0,t_1,\cdots,t_k}={\rm Exp}(\xi_0)\cdot\gamma_{\xi^1(t_1)}\cdot\,\cdots\,\cdot\gamma_{\xi^k(t_k)}$ 
belongs to $\Sigma$.  
\hspace{8truecm}q.e.d.

\vspace{0.5truecm}

It is clear that any point of $\Sigma$ is described as $x_{\xi_0,t_1,\cdots,t_k}$ for some 
$\xi_0\in\mathfrak b$ and some $t_1,\cdots,t_k\in\mathbb R$.  
Fix an orthonormal base $\{e^0_1,\cdots,e^0_{m_0}\}$ of $\mathfrak b$, where $m_0:={\rm dim}\,\mathfrak b$.  
Define vector fields $E^0_i$ ($i=1,\cdots,m_0$) and $E^j$ ($j=1,\cdots,k$) along $\Sigma$ by 
$$\begin{array}{l}
\hspace{1.5truecm}\displaystyle{(E^0_i)_{x_{\xi_0,t_1,\cdots,t_k}}:=(L_{x_{\xi_0,t_1,\cdots,t_k}})_{\ast}(e^0_i)
(=(\widetilde{e^0_i})_{x_{\xi_0,t_1,\cdots,t_k}})}\\
{\rm and}\qquad\,\,
\displaystyle{(E^j)_{x_{\xi_0,t_1,\cdots,t_k}}:=(L_{x_{\xi_0,t_1,\cdots,t_k}})_{\ast}(\xi^j_{t_j})
(=(\widetilde{\xi^j_{t_j}})_{x_{\xi_0,t_1,\cdots,t_k}}).}
\end{array}$$
By imitating the discussions in the proofs of Lemmas 4.1 and 4.2, 
we can show the following fact for these vector fields.  

\vspace{0.5truecm}

\noindent
{\bf Lemma 4.3.} {\sl The vector fields $E^0_i$ ($i=1,\cdots,m_0$) and $E^j$ ($j=1,\cdots,k$) are tangent to 
$\Sigma$ and they give a parallel orthonormal tangent frame field on $\Sigma$.}

\vspace{0.5truecm}

\noindent
{\it Proof.} Let $(\widehat{\xi^i})^j$ (resp. $(\widehat{\xi^i})^0$) be the parallel vector field along 
$\gamma_{\xi^j}$ ($i\not=j$) (resp. $\gamma_{\xi_0}$) with $(\widehat{\xi^i})^j_0=\xi^i$ 
(resp. $(\widehat{\xi^i})^0_0=\xi^i$) and 
$(\widehat{\xi_0})^j$ be the parallel vector field along 
$\gamma_{\xi^j}$ with $(\widehat{\xi_0})^j_0=\xi_0$.  
According to Lemma 4.1, we have $(\gamma_{\xi^i})'(t)=(L_{\gamma_{\xi^i}(t)})_{\ast}(\xi^i_t)$ and 
$(\gamma_{\xi_0})'(t)=(L_{\gamma_{\xi_0}(t)})_{\ast}(\xi_0)$.  
Also, we can show 
$(\widehat{\xi^i})^j_{\gamma_{\xi^j}(t)}=(L_{\gamma_{\xi^j}(t)})_{\ast}(\xi^i)$ ($j\not=i$), 
$(\widehat{\xi^i})^0_{\gamma_{\xi_0}(t)}=(L_{\gamma_{\xi_0}(t)})_{\ast}(\xi^i)$ and 
$(\widehat{\xi_0})^j_{\gamma_{\xi^j}(t)}=(L_{\gamma_{\xi^j}(t)})_{\ast}(\xi_0)$ 
by imitating the discussion in the proof of Lemma 4.2.  
On the basis of these facts, we can derive the statement of this lemma, where we note that $\Sigma$ is flat.  
\hspace{12.8truecm}q.e.d.

\vspace{0.5truecm}

\centerline{
\unitlength 0.1in
\begin{picture}( 53.7700, 26.0500)(  1.5000,-33.8500)
%
\special{pn 8}%
\special{pa 1610 1190}%
\special{pa 1048 1756}%
\special{pa 2820 1756}%
\special{pa 3310 1190}%
\special{pa 3310 1190}%
\special{pa 1610 1190}%
\special{fp}%
%
\special{pn 20}%
\special{sh 1}%
\special{ar 1530 1604 10 10 0  6.28318530717959E+0000}%
\special{sh 1}%
\special{ar 1530 1604 10 10 0  6.28318530717959E+0000}%
%
\special{pn 13}%
\special{pa 1540 1604}%
\special{pa 1816 1604}%
\special{fp}%
\special{sh 1}%
\special{pa 1816 1604}%
\special{pa 1748 1584}%
\special{pa 1762 1604}%
\special{pa 1748 1624}%
\special{pa 1816 1604}%
\special{fp}%
%
\special{pn 8}%
\special{pa 1216 1596}%
\special{pa 2958 1596}%
\special{fp}%
%
\special{pn 8}%
\special{pa 1510 1294}%
\special{pa 3232 1294}%
\special{fp}%
%
\special{pn 8}%
\special{pa 1364 1756}%
\special{pa 1972 1190}%
\special{fp}%
%
\special{pn 13}%
\special{pa 2440 1606}%
\special{pa 2626 1798}%
\special{fp}%
\special{sh 1}%
\special{pa 2626 1798}%
\special{pa 2594 1736}%
\special{pa 2590 1760}%
\special{pa 2566 1764}%
\special{pa 2626 1798}%
\special{fp}%
%
\special{pn 13}%
\special{pa 1856 1304}%
\special{pa 2130 1304}%
\special{fp}%
\special{sh 1}%
\special{pa 2130 1304}%
\special{pa 2064 1284}%
\special{pa 2078 1304}%
\special{pa 2064 1324}%
\special{pa 2130 1304}%
\special{fp}%
%
\special{pn 8}%
\special{pa 2258 1756}%
\special{pa 2868 1190}%
\special{fp}%
%
\special{pn 8}%
\special{ar 1530 1604 296 282  0.0000000 0.0415225}%
\special{ar 1530 1604 296 282  0.1660900 0.2076125}%
\special{ar 1530 1604 296 282  0.3321799 0.3737024}%
\special{ar 1530 1604 296 282  0.4982699 0.5397924}%
\special{ar 1530 1604 296 282  0.6643599 0.7058824}%
\special{ar 1530 1604 296 282  0.8304498 0.8719723}%
\special{ar 1530 1604 296 282  0.9965398 1.0380623}%
\special{ar 1530 1604 296 282  1.1626298 1.2041522}%
\special{ar 1530 1604 296 282  1.3287197 1.3702422}%
\special{ar 1530 1604 296 282  1.4948097 1.5363322}%
\special{ar 1530 1604 296 282  1.6608997 1.7024221}%
\special{ar 1530 1604 296 282  1.8269896 1.8685121}%
\special{ar 1530 1604 296 282  1.9930796 2.0346021}%
\special{ar 1530 1604 296 282  2.1591696 2.2006920}%
\special{ar 1530 1604 296 282  2.3252595 2.3667820}%
\special{ar 1530 1604 296 282  2.4913495 2.5328720}%
\special{ar 1530 1604 296 282  2.6574394 2.6989619}%
\special{ar 1530 1604 296 282  2.8235294 2.8650519}%
\special{ar 1530 1604 296 282  2.9896194 3.0311419}%
\special{ar 1530 1604 296 282  3.1557093 3.1972318}%
\special{ar 1530 1604 296 282  3.3217993 3.3633218}%
\special{ar 1530 1604 296 282  3.4878893 3.5294118}%
\special{ar 1530 1604 296 282  3.6539792 3.6955017}%
\special{ar 1530 1604 296 282  3.8200692 3.8615917}%
\special{ar 1530 1604 296 282  3.9861592 4.0276817}%
\special{ar 1530 1604 296 282  4.1522491 4.1937716}%
\special{ar 1530 1604 296 282  4.3183391 4.3598616}%
\special{ar 1530 1604 296 282  4.4844291 4.5259516}%
\special{ar 1530 1604 296 282  4.6505190 4.6920415}%
\special{ar 1530 1604 296 282  4.8166090 4.8581315}%
\special{ar 1530 1604 296 282  4.9826990 5.0242215}%
\special{ar 1530 1604 296 282  5.1487889 5.1903114}%
\special{ar 1530 1604 296 282  5.3148789 5.3564014}%
\special{ar 1530 1604 296 282  5.4809689 5.5224913}%
\special{ar 1530 1604 296 282  5.6470588 5.6885813}%
\special{ar 1530 1604 296 282  5.8131488 5.8546713}%
\special{ar 1530 1604 296 282  5.9792388 6.0207612}%
\special{ar 1530 1604 296 282  6.1453287 6.1868512}%
%
\special{pn 8}%
\special{ar 1846 1304 296 282  0.0000000 0.0415945}%
\special{ar 1846 1304 296 282  0.1663778 0.2079723}%
\special{ar 1846 1304 296 282  0.3327556 0.3743501}%
\special{ar 1846 1304 296 282  0.4991334 0.5407279}%
\special{ar 1846 1304 296 282  0.6655113 0.7071057}%
\special{ar 1846 1304 296 282  0.8318891 0.8734835}%
\special{ar 1846 1304 296 282  0.9982669 1.0398614}%
\special{ar 1846 1304 296 282  1.1646447 1.2062392}%
\special{ar 1846 1304 296 282  1.3310225 1.3726170}%
\special{ar 1846 1304 296 282  1.4974003 1.5389948}%
\special{ar 1846 1304 296 282  1.6637782 1.7053726}%
\special{ar 1846 1304 296 282  1.8301560 1.8717504}%
\special{ar 1846 1304 296 282  1.9965338 2.0381282}%
\special{ar 1846 1304 296 282  2.1629116 2.2045061}%
\special{ar 1846 1304 296 282  2.3292894 2.3708839}%
\special{ar 1846 1304 296 282  2.4956672 2.5372617}%
\special{ar 1846 1304 296 282  2.6620451 2.7036395}%
\special{ar 1846 1304 296 282  2.8284229 2.8700173}%
\special{ar 1846 1304 296 282  2.9948007 3.0363951}%
\special{ar 1846 1304 296 282  3.1611785 3.2027730}%
\special{ar 1846 1304 296 282  3.3275563 3.3691508}%
\special{ar 1846 1304 296 282  3.4939341 3.5355286}%
\special{ar 1846 1304 296 282  3.6603120 3.7019064}%
\special{ar 1846 1304 296 282  3.8266898 3.8682842}%
\special{ar 1846 1304 296 282  3.9930676 4.0346620}%
\special{ar 1846 1304 296 282  4.1594454 4.2010399}%
\special{ar 1846 1304 296 282  4.3258232 4.3674177}%
\special{ar 1846 1304 296 282  4.4922010 4.5337955}%
\special{ar 1846 1304 296 282  4.6585789 4.7001733}%
\special{ar 1846 1304 296 282  4.8249567 4.8665511}%
\special{ar 1846 1304 296 282  4.9913345 5.0329289}%
\special{ar 1846 1304 296 282  5.1577123 5.1993068}%
\special{ar 1846 1304 296 282  5.3240901 5.3656846}%
\special{ar 1846 1304 296 282  5.4904679 5.5320624}%
\special{ar 1846 1304 296 282  5.6568458 5.6984402}%
\special{ar 1846 1304 296 282  5.8232236 5.8648180}%
\special{ar 1846 1304 296 282  5.9896014 6.0311958}%
\special{ar 1846 1304 296 282  6.1559792 6.1975737}%
%
\special{pn 8}%
\special{pa 3826 1200}%
\special{pa 3264 1764}%
\special{pa 5036 1764}%
\special{pa 5528 1200}%
\special{pa 5528 1200}%
\special{pa 3826 1200}%
\special{fp}%
%
\special{pn 20}%
\special{sh 1}%
\special{ar 3746 1614 10 10 0  6.28318530717959E+0000}%
\special{sh 1}%
\special{ar 3746 1614 10 10 0  6.28318530717959E+0000}%
%
\special{pn 8}%
\special{pa 3422 1614}%
\special{pa 5164 1614}%
\special{fp}%
%
\special{pn 8}%
\special{pa 3708 1312}%
\special{pa 5448 1312}%
\special{fp}%
%
\special{pn 8}%
\special{pa 3580 1764}%
\special{pa 4170 1180}%
\special{fp}%
%
\special{pn 13}%
\special{pa 4652 1614}%
\special{pa 4928 1614}%
\special{fp}%
\special{sh 1}%
\special{pa 4928 1614}%
\special{pa 4860 1594}%
\special{pa 4874 1614}%
\special{pa 4860 1634}%
\special{pa 4928 1614}%
\special{fp}%
%
\special{pn 13}%
\special{pa 4976 1312}%
\special{pa 5252 1312}%
\special{fp}%
\special{sh 1}%
\special{pa 5252 1312}%
\special{pa 5186 1292}%
\special{pa 5200 1312}%
\special{pa 5186 1332}%
\special{pa 5252 1312}%
\special{fp}%
%
\special{pn 8}%
\special{pa 4474 1764}%
\special{pa 5084 1200}%
\special{fp}%
%
\special{pn 8}%
\special{ar 3746 1614 296 282  0.0000000 0.0415225}%
\special{ar 3746 1614 296 282  0.1660900 0.2076125}%
\special{ar 3746 1614 296 282  0.3321799 0.3737024}%
\special{ar 3746 1614 296 282  0.4982699 0.5397924}%
\special{ar 3746 1614 296 282  0.6643599 0.7058824}%
\special{ar 3746 1614 296 282  0.8304498 0.8719723}%
\special{ar 3746 1614 296 282  0.9965398 1.0380623}%
\special{ar 3746 1614 296 282  1.1626298 1.2041522}%
\special{ar 3746 1614 296 282  1.3287197 1.3702422}%
\special{ar 3746 1614 296 282  1.4948097 1.5363322}%
\special{ar 3746 1614 296 282  1.6608997 1.7024221}%
\special{ar 3746 1614 296 282  1.8269896 1.8685121}%
\special{ar 3746 1614 296 282  1.9930796 2.0346021}%
\special{ar 3746 1614 296 282  2.1591696 2.2006920}%
\special{ar 3746 1614 296 282  2.3252595 2.3667820}%
\special{ar 3746 1614 296 282  2.4913495 2.5328720}%
\special{ar 3746 1614 296 282  2.6574394 2.6989619}%
\special{ar 3746 1614 296 282  2.8235294 2.8650519}%
\special{ar 3746 1614 296 282  2.9896194 3.0311419}%
\special{ar 3746 1614 296 282  3.1557093 3.1972318}%
\special{ar 3746 1614 296 282  3.3217993 3.3633218}%
\special{ar 3746 1614 296 282  3.4878893 3.5294118}%
\special{ar 3746 1614 296 282  3.6539792 3.6955017}%
\special{ar 3746 1614 296 282  3.8200692 3.8615917}%
\special{ar 3746 1614 296 282  3.9861592 4.0276817}%
\special{ar 3746 1614 296 282  4.1522491 4.1937716}%
\special{ar 3746 1614 296 282  4.3183391 4.3598616}%
\special{ar 3746 1614 296 282  4.4844291 4.5259516}%
\special{ar 3746 1614 296 282  4.6505190 4.6920415}%
\special{ar 3746 1614 296 282  4.8166090 4.8581315}%
\special{ar 3746 1614 296 282  4.9826990 5.0242215}%
\special{ar 3746 1614 296 282  5.1487889 5.1903114}%
\special{ar 3746 1614 296 282  5.3148789 5.3564014}%
\special{ar 3746 1614 296 282  5.4809689 5.5224913}%
\special{ar 3746 1614 296 282  5.6470588 5.6885813}%
\special{ar 3746 1614 296 282  5.8131488 5.8546713}%
\special{ar 3746 1614 296 282  5.9792388 6.0207612}%
\special{ar 3746 1614 296 282  6.1453287 6.1868512}%
%
\special{pn 8}%
\special{ar 4042 1312 296 284  0.0000000 0.0415225}%
\special{ar 4042 1312 296 284  0.1660900 0.2076125}%
\special{ar 4042 1312 296 284  0.3321799 0.3737024}%
\special{ar 4042 1312 296 284  0.4982699 0.5397924}%
\special{ar 4042 1312 296 284  0.6643599 0.7058824}%
\special{ar 4042 1312 296 284  0.8304498 0.8719723}%
\special{ar 4042 1312 296 284  0.9965398 1.0380623}%
\special{ar 4042 1312 296 284  1.1626298 1.2041522}%
\special{ar 4042 1312 296 284  1.3287197 1.3702422}%
\special{ar 4042 1312 296 284  1.4948097 1.5363322}%
\special{ar 4042 1312 296 284  1.6608997 1.7024221}%
\special{ar 4042 1312 296 284  1.8269896 1.8685121}%
\special{ar 4042 1312 296 284  1.9930796 2.0346021}%
\special{ar 4042 1312 296 284  2.1591696 2.2006920}%
\special{ar 4042 1312 296 284  2.3252595 2.3667820}%
\special{ar 4042 1312 296 284  2.4913495 2.5328720}%
\special{ar 4042 1312 296 284  2.6574394 2.6989619}%
\special{ar 4042 1312 296 284  2.8235294 2.8650519}%
\special{ar 4042 1312 296 284  2.9896194 3.0311419}%
\special{ar 4042 1312 296 284  3.1557093 3.1972318}%
\special{ar 4042 1312 296 284  3.3217993 3.3633218}%
\special{ar 4042 1312 296 284  3.4878893 3.5294118}%
\special{ar 4042 1312 296 284  3.6539792 3.6955017}%
\special{ar 4042 1312 296 284  3.8200692 3.8615917}%
\special{ar 4042 1312 296 284  3.9861592 4.0276817}%
\special{ar 4042 1312 296 284  4.1522491 4.1937716}%
\special{ar 4042 1312 296 284  4.3183391 4.3598616}%
\special{ar 4042 1312 296 284  4.4844291 4.5259516}%
\special{ar 4042 1312 296 284  4.6505190 4.6920415}%
\special{ar 4042 1312 296 284  4.8166090 4.8581315}%
\special{ar 4042 1312 296 284  4.9826990 5.0242215}%
\special{ar 4042 1312 296 284  5.1487889 5.1903114}%
\special{ar 4042 1312 296 284  5.3148789 5.3564014}%
\special{ar 4042 1312 296 284  5.4809689 5.5224913}%
\special{ar 4042 1312 296 284  5.6470588 5.6885813}%
\special{ar 4042 1312 296 284  5.8131488 5.8546713}%
\special{ar 4042 1312 296 284  5.9792388 6.0207612}%
\special{ar 4042 1312 296 284  6.1453287 6.1868512}%
%
\special{pn 13}%
\special{pa 3746 1614}%
\special{pa 3978 1446}%
\special{fp}%
\special{sh 1}%
\special{pa 3978 1446}%
\special{pa 3912 1468}%
\special{pa 3936 1476}%
\special{pa 3936 1500}%
\special{pa 3978 1446}%
\special{fp}%
%
\special{pn 13}%
\special{pa 4042 1312}%
\special{pa 4268 1138}%
\special{fp}%
\special{sh 1}%
\special{pa 4268 1138}%
\special{pa 4204 1164}%
\special{pa 4226 1170}%
\special{pa 4228 1194}%
\special{pa 4268 1138}%
\special{fp}%
%
\special{pn 8}%
\special{ar 2426 1604 296 282  0.0000000 0.0415945}%
\special{ar 2426 1604 296 282  0.1663778 0.2079723}%
\special{ar 2426 1604 296 282  0.3327556 0.3743501}%
\special{ar 2426 1604 296 282  0.4991334 0.5407279}%
\special{ar 2426 1604 296 282  0.6655113 0.7071057}%
\special{ar 2426 1604 296 282  0.8318891 0.8734835}%
\special{ar 2426 1604 296 282  0.9982669 1.0398614}%
\special{ar 2426 1604 296 282  1.1646447 1.2062392}%
\special{ar 2426 1604 296 282  1.3310225 1.3726170}%
\special{ar 2426 1604 296 282  1.4974003 1.5389948}%
\special{ar 2426 1604 296 282  1.6637782 1.7053726}%
\special{ar 2426 1604 296 282  1.8301560 1.8717504}%
\special{ar 2426 1604 296 282  1.9965338 2.0381282}%
\special{ar 2426 1604 296 282  2.1629116 2.2045061}%
\special{ar 2426 1604 296 282  2.3292894 2.3708839}%
\special{ar 2426 1604 296 282  2.4956672 2.5372617}%
\special{ar 2426 1604 296 282  2.6620451 2.7036395}%
\special{ar 2426 1604 296 282  2.8284229 2.8700173}%
\special{ar 2426 1604 296 282  2.9948007 3.0363951}%
\special{ar 2426 1604 296 282  3.1611785 3.2027730}%
\special{ar 2426 1604 296 282  3.3275563 3.3691508}%
\special{ar 2426 1604 296 282  3.4939341 3.5355286}%
\special{ar 2426 1604 296 282  3.6603120 3.7019064}%
\special{ar 2426 1604 296 282  3.8266898 3.8682842}%
\special{ar 2426 1604 296 282  3.9930676 4.0346620}%
\special{ar 2426 1604 296 282  4.1594454 4.2010399}%
\special{ar 2426 1604 296 282  4.3258232 4.3674177}%
\special{ar 2426 1604 296 282  4.4922010 4.5337955}%
\special{ar 2426 1604 296 282  4.6585789 4.7001733}%
\special{ar 2426 1604 296 282  4.8249567 4.8665511}%
\special{ar 2426 1604 296 282  4.9913345 5.0329289}%
\special{ar 2426 1604 296 282  5.1577123 5.1993068}%
\special{ar 2426 1604 296 282  5.3240901 5.3656846}%
\special{ar 2426 1604 296 282  5.4904679 5.5320624}%
\special{ar 2426 1604 296 282  5.6568458 5.6984402}%
\special{ar 2426 1604 296 282  5.8232236 5.8648180}%
\special{ar 2426 1604 296 282  5.9896014 6.0311958}%
\special{ar 2426 1604 296 282  6.1559792 6.1975737}%
%
\special{pn 8}%
\special{ar 2752 1304 294 282  0.0000000 0.0416667}%
\special{ar 2752 1304 294 282  0.1666667 0.2083333}%
\special{ar 2752 1304 294 282  0.3333333 0.3750000}%
\special{ar 2752 1304 294 282  0.5000000 0.5416667}%
\special{ar 2752 1304 294 282  0.6666667 0.7083333}%
\special{ar 2752 1304 294 282  0.8333333 0.8750000}%
\special{ar 2752 1304 294 282  1.0000000 1.0416667}%
\special{ar 2752 1304 294 282  1.1666667 1.2083333}%
\special{ar 2752 1304 294 282  1.3333333 1.3750000}%
\special{ar 2752 1304 294 282  1.5000000 1.5416667}%
\special{ar 2752 1304 294 282  1.6666667 1.7083333}%
\special{ar 2752 1304 294 282  1.8333333 1.8750000}%
\special{ar 2752 1304 294 282  2.0000000 2.0416667}%
\special{ar 2752 1304 294 282  2.1666667 2.2083333}%
\special{ar 2752 1304 294 282  2.3333333 2.3750000}%
\special{ar 2752 1304 294 282  2.5000000 2.5416667}%
\special{ar 2752 1304 294 282  2.6666667 2.7083333}%
\special{ar 2752 1304 294 282  2.8333333 2.8750000}%
\special{ar 2752 1304 294 282  3.0000000 3.0416667}%
\special{ar 2752 1304 294 282  3.1666667 3.2083333}%
\special{ar 2752 1304 294 282  3.3333333 3.3750000}%
\special{ar 2752 1304 294 282  3.5000000 3.5416667}%
\special{ar 2752 1304 294 282  3.6666667 3.7083333}%
\special{ar 2752 1304 294 282  3.8333333 3.8750000}%
\special{ar 2752 1304 294 282  4.0000000 4.0416667}%
\special{ar 2752 1304 294 282  4.1666667 4.2083333}%
\special{ar 2752 1304 294 282  4.3333333 4.3750000}%
\special{ar 2752 1304 294 282  4.5000000 4.5416667}%
\special{ar 2752 1304 294 282  4.6666667 4.7083333}%
\special{ar 2752 1304 294 282  4.8333333 4.8750000}%
\special{ar 2752 1304 294 282  5.0000000 5.0416667}%
\special{ar 2752 1304 294 282  5.1666667 5.2083333}%
\special{ar 2752 1304 294 282  5.3333333 5.3750000}%
\special{ar 2752 1304 294 282  5.5000000 5.5416667}%
\special{ar 2752 1304 294 282  5.6666667 5.7083333}%
\special{ar 2752 1304 294 282  5.8333333 5.8750000}%
\special{ar 2752 1304 294 282  6.0000000 6.0416667}%
\special{ar 2752 1304 294 282  6.1666667 6.2083333}%
%
\special{pn 13}%
\special{pa 2760 1304}%
\special{pa 2948 1496}%
\special{fp}%
\special{sh 1}%
\special{pa 2948 1496}%
\special{pa 2916 1434}%
\special{pa 2910 1458}%
\special{pa 2886 1462}%
\special{pa 2948 1496}%
\special{fp}%
%
\special{pn 8}%
\special{ar 4642 1614 294 282  0.0000000 0.0416667}%
\special{ar 4642 1614 294 282  0.1666667 0.2083333}%
\special{ar 4642 1614 294 282  0.3333333 0.3750000}%
\special{ar 4642 1614 294 282  0.5000000 0.5416667}%
\special{ar 4642 1614 294 282  0.6666667 0.7083333}%
\special{ar 4642 1614 294 282  0.8333333 0.8750000}%
\special{ar 4642 1614 294 282  1.0000000 1.0416667}%
\special{ar 4642 1614 294 282  1.1666667 1.2083333}%
\special{ar 4642 1614 294 282  1.3333333 1.3750000}%
\special{ar 4642 1614 294 282  1.5000000 1.5416667}%
\special{ar 4642 1614 294 282  1.6666667 1.7083333}%
\special{ar 4642 1614 294 282  1.8333333 1.8750000}%
\special{ar 4642 1614 294 282  2.0000000 2.0416667}%
\special{ar 4642 1614 294 282  2.1666667 2.2083333}%
\special{ar 4642 1614 294 282  2.3333333 2.3750000}%
\special{ar 4642 1614 294 282  2.5000000 2.5416667}%
\special{ar 4642 1614 294 282  2.6666667 2.7083333}%
\special{ar 4642 1614 294 282  2.8333333 2.8750000}%
\special{ar 4642 1614 294 282  3.0000000 3.0416667}%
\special{ar 4642 1614 294 282  3.1666667 3.2083333}%
\special{ar 4642 1614 294 282  3.3333333 3.3750000}%
\special{ar 4642 1614 294 282  3.5000000 3.5416667}%
\special{ar 4642 1614 294 282  3.6666667 3.7083333}%
\special{ar 4642 1614 294 282  3.8333333 3.8750000}%
\special{ar 4642 1614 294 282  4.0000000 4.0416667}%
\special{ar 4642 1614 294 282  4.1666667 4.2083333}%
\special{ar 4642 1614 294 282  4.3333333 4.3750000}%
\special{ar 4642 1614 294 282  4.5000000 4.5416667}%
\special{ar 4642 1614 294 282  4.6666667 4.7083333}%
\special{ar 4642 1614 294 282  4.8333333 4.8750000}%
\special{ar 4642 1614 294 282  5.0000000 5.0416667}%
\special{ar 4642 1614 294 282  5.1666667 5.2083333}%
\special{ar 4642 1614 294 282  5.3333333 5.3750000}%
\special{ar 4642 1614 294 282  5.5000000 5.5416667}%
\special{ar 4642 1614 294 282  5.6666667 5.7083333}%
\special{ar 4642 1614 294 282  5.8333333 5.8750000}%
\special{ar 4642 1614 294 282  6.0000000 6.0416667}%
\special{ar 4642 1614 294 282  6.1666667 6.2083333}%
%
\special{pn 8}%
\special{ar 4968 1312 294 284  0.0000000 0.0415945}%
\special{ar 4968 1312 294 284  0.1663778 0.2079723}%
\special{ar 4968 1312 294 284  0.3327556 0.3743501}%
\special{ar 4968 1312 294 284  0.4991334 0.5407279}%
\special{ar 4968 1312 294 284  0.6655113 0.7071057}%
\special{ar 4968 1312 294 284  0.8318891 0.8734835}%
\special{ar 4968 1312 294 284  0.9982669 1.0398614}%
\special{ar 4968 1312 294 284  1.1646447 1.2062392}%
\special{ar 4968 1312 294 284  1.3310225 1.3726170}%
\special{ar 4968 1312 294 284  1.4974003 1.5389948}%
\special{ar 4968 1312 294 284  1.6637782 1.7053726}%
\special{ar 4968 1312 294 284  1.8301560 1.8717504}%
\special{ar 4968 1312 294 284  1.9965338 2.0381282}%
\special{ar 4968 1312 294 284  2.1629116 2.2045061}%
\special{ar 4968 1312 294 284  2.3292894 2.3708839}%
\special{ar 4968 1312 294 284  2.4956672 2.5372617}%
\special{ar 4968 1312 294 284  2.6620451 2.7036395}%
\special{ar 4968 1312 294 284  2.8284229 2.8700173}%
\special{ar 4968 1312 294 284  2.9948007 3.0363951}%
\special{ar 4968 1312 294 284  3.1611785 3.2027730}%
\special{ar 4968 1312 294 284  3.3275563 3.3691508}%
\special{ar 4968 1312 294 284  3.4939341 3.5355286}%
\special{ar 4968 1312 294 284  3.6603120 3.7019064}%
\special{ar 4968 1312 294 284  3.8266898 3.8682842}%
\special{ar 4968 1312 294 284  3.9930676 4.0346620}%
\special{ar 4968 1312 294 284  4.1594454 4.2010399}%
\special{ar 4968 1312 294 284  4.3258232 4.3674177}%
\special{ar 4968 1312 294 284  4.4922010 4.5337955}%
\special{ar 4968 1312 294 284  4.6585789 4.7001733}%
\special{ar 4968 1312 294 284  4.8249567 4.8665511}%
\special{ar 4968 1312 294 284  4.9913345 5.0329289}%
\special{ar 4968 1312 294 284  5.1577123 5.1993068}%
\special{ar 4968 1312 294 284  5.3240901 5.3656846}%
\special{ar 4968 1312 294 284  5.4904679 5.5320624}%
\special{ar 4968 1312 294 284  5.6568458 5.6984402}%
\special{ar 4968 1312 294 284  5.8232236 5.8648180}%
\special{ar 4968 1312 294 284  5.9896014 6.0311958}%
\special{ar 4968 1312 294 284  6.1559792 6.1975737}%
%
\special{pn 8}%
\special{pa 2460 2576}%
\special{pa 1882 3188}%
\special{pa 3710 3188}%
\special{pa 4218 2576}%
\special{pa 4218 2576}%
\special{pa 2460 2576}%
\special{fp}%
%
\special{pn 20}%
\special{sh 1}%
\special{ar 2380 3026 10 10 0  6.28318530717959E+0000}%
\special{sh 1}%
\special{ar 2380 3026 10 10 0  6.28318530717959E+0000}%
%
\special{pn 13}%
\special{pa 2388 3026}%
\special{pa 2674 3026}%
\special{fp}%
\special{sh 1}%
\special{pa 2674 3026}%
\special{pa 2608 3006}%
\special{pa 2622 3026}%
\special{pa 2608 3046}%
\special{pa 2674 3026}%
\special{fp}%
%
\special{pn 8}%
\special{pa 2054 3016}%
\special{pa 3852 3016}%
\special{fp}%
%
\special{pn 8}%
\special{pa 2358 2688}%
\special{pa 4136 2688}%
\special{fp}%
%
\special{pn 8}%
\special{pa 2206 3188}%
\special{pa 2836 2576}%
\special{fp}%
%
\special{pn 13}%
\special{pa 2714 2698}%
\special{pa 3000 2698}%
\special{fp}%
\special{sh 1}%
\special{pa 3000 2698}%
\special{pa 2932 2678}%
\special{pa 2946 2698}%
\special{pa 2932 2718}%
\special{pa 3000 2698}%
\special{fp}%
%
\special{pn 8}%
\special{pa 3130 3188}%
\special{pa 3760 2576}%
\special{fp}%
%
\special{pn 8}%
\special{ar 2380 3026 304 306  0.0000000 0.0393443}%
\special{ar 2380 3026 304 306  0.1573770 0.1967213}%
\special{ar 2380 3026 304 306  0.3147541 0.3540984}%
\special{ar 2380 3026 304 306  0.4721311 0.5114754}%
\special{ar 2380 3026 304 306  0.6295082 0.6688525}%
\special{ar 2380 3026 304 306  0.7868852 0.8262295}%
\special{ar 2380 3026 304 306  0.9442623 0.9836066}%
\special{ar 2380 3026 304 306  1.1016393 1.1409836}%
\special{ar 2380 3026 304 306  1.2590164 1.2983607}%
\special{ar 2380 3026 304 306  1.4163934 1.4557377}%
\special{ar 2380 3026 304 306  1.5737705 1.6131148}%
\special{ar 2380 3026 304 306  1.7311475 1.7704918}%
\special{ar 2380 3026 304 306  1.8885246 1.9278689}%
\special{ar 2380 3026 304 306  2.0459016 2.0852459}%
\special{ar 2380 3026 304 306  2.2032787 2.2426230}%
\special{ar 2380 3026 304 306  2.3606557 2.4000000}%
\special{ar 2380 3026 304 306  2.5180328 2.5573770}%
\special{ar 2380 3026 304 306  2.6754098 2.7147541}%
\special{ar 2380 3026 304 306  2.8327869 2.8721311}%
\special{ar 2380 3026 304 306  2.9901639 3.0295082}%
\special{ar 2380 3026 304 306  3.1475410 3.1868852}%
\special{ar 2380 3026 304 306  3.3049180 3.3442623}%
\special{ar 2380 3026 304 306  3.4622951 3.5016393}%
\special{ar 2380 3026 304 306  3.6196721 3.6590164}%
\special{ar 2380 3026 304 306  3.7770492 3.8163934}%
\special{ar 2380 3026 304 306  3.9344262 3.9737705}%
\special{ar 2380 3026 304 306  4.0918033 4.1311475}%
\special{ar 2380 3026 304 306  4.2491803 4.2885246}%
\special{ar 2380 3026 304 306  4.4065574 4.4459016}%
\special{ar 2380 3026 304 306  4.5639344 4.6032787}%
\special{ar 2380 3026 304 306  4.7213115 4.7606557}%
\special{ar 2380 3026 304 306  4.8786885 4.9180328}%
\special{ar 2380 3026 304 306  5.0360656 5.0754098}%
\special{ar 2380 3026 304 306  5.1934426 5.2327869}%
\special{ar 2380 3026 304 306  5.3508197 5.3901639}%
\special{ar 2380 3026 304 306  5.5081967 5.5475410}%
\special{ar 2380 3026 304 306  5.6655738 5.7049180}%
\special{ar 2380 3026 304 306  5.8229508 5.8622951}%
\special{ar 2380 3026 304 306  5.9803279 6.0196721}%
\special{ar 2380 3026 304 306  6.1377049 6.1770492}%
%
\special{pn 8}%
\special{ar 2704 2698 304 306  0.0000000 0.0393443}%
\special{ar 2704 2698 304 306  0.1573770 0.1967213}%
\special{ar 2704 2698 304 306  0.3147541 0.3540984}%
\special{ar 2704 2698 304 306  0.4721311 0.5114754}%
\special{ar 2704 2698 304 306  0.6295082 0.6688525}%
\special{ar 2704 2698 304 306  0.7868852 0.8262295}%
\special{ar 2704 2698 304 306  0.9442623 0.9836066}%
\special{ar 2704 2698 304 306  1.1016393 1.1409836}%
\special{ar 2704 2698 304 306  1.2590164 1.2983607}%
\special{ar 2704 2698 304 306  1.4163934 1.4557377}%
\special{ar 2704 2698 304 306  1.5737705 1.6131148}%
\special{ar 2704 2698 304 306  1.7311475 1.7704918}%
\special{ar 2704 2698 304 306  1.8885246 1.9278689}%
\special{ar 2704 2698 304 306  2.0459016 2.0852459}%
\special{ar 2704 2698 304 306  2.2032787 2.2426230}%
\special{ar 2704 2698 304 306  2.3606557 2.4000000}%
\special{ar 2704 2698 304 306  2.5180328 2.5573770}%
\special{ar 2704 2698 304 306  2.6754098 2.7147541}%
\special{ar 2704 2698 304 306  2.8327869 2.8721311}%
\special{ar 2704 2698 304 306  2.9901639 3.0295082}%
\special{ar 2704 2698 304 306  3.1475410 3.1868852}%
\special{ar 2704 2698 304 306  3.3049180 3.3442623}%
\special{ar 2704 2698 304 306  3.4622951 3.5016393}%
\special{ar 2704 2698 304 306  3.6196721 3.6590164}%
\special{ar 2704 2698 304 306  3.7770492 3.8163934}%
\special{ar 2704 2698 304 306  3.9344262 3.9737705}%
\special{ar 2704 2698 304 306  4.0918033 4.1311475}%
\special{ar 2704 2698 304 306  4.2491803 4.2885246}%
\special{ar 2704 2698 304 306  4.4065574 4.4459016}%
\special{ar 2704 2698 304 306  4.5639344 4.6032787}%
\special{ar 2704 2698 304 306  4.7213115 4.7606557}%
\special{ar 2704 2698 304 306  4.8786885 4.9180328}%
\special{ar 2704 2698 304 306  5.0360656 5.0754098}%
\special{ar 2704 2698 304 306  5.1934426 5.2327869}%
\special{ar 2704 2698 304 306  5.3508197 5.3901639}%
\special{ar 2704 2698 304 306  5.5081967 5.5475410}%
\special{ar 2704 2698 304 306  5.6655738 5.7049180}%
\special{ar 2704 2698 304 306  5.8229508 5.8622951}%
\special{ar 2704 2698 304 306  5.9803279 6.0196721}%
\special{ar 2704 2698 304 306  6.1377049 6.1770492}%
%
\special{pn 8}%
\special{ar 3304 3026 306 306  0.0000000 0.0392799}%
\special{ar 3304 3026 306 306  0.1571195 0.1963993}%
\special{ar 3304 3026 306 306  0.3142390 0.3535188}%
\special{ar 3304 3026 306 306  0.4713584 0.5106383}%
\special{ar 3304 3026 306 306  0.6284779 0.6677578}%
\special{ar 3304 3026 306 306  0.7855974 0.8248773}%
\special{ar 3304 3026 306 306  0.9427169 0.9819967}%
\special{ar 3304 3026 306 306  1.0998363 1.1391162}%
\special{ar 3304 3026 306 306  1.2569558 1.2962357}%
\special{ar 3304 3026 306 306  1.4140753 1.4533552}%
\special{ar 3304 3026 306 306  1.5711948 1.6104746}%
\special{ar 3304 3026 306 306  1.7283142 1.7675941}%
\special{ar 3304 3026 306 306  1.8854337 1.9247136}%
\special{ar 3304 3026 306 306  2.0425532 2.0818331}%
\special{ar 3304 3026 306 306  2.1996727 2.2389525}%
\special{ar 3304 3026 306 306  2.3567921 2.3960720}%
\special{ar 3304 3026 306 306  2.5139116 2.5531915}%
\special{ar 3304 3026 306 306  2.6710311 2.7103110}%
\special{ar 3304 3026 306 306  2.8281506 2.8674304}%
\special{ar 3304 3026 306 306  2.9852700 3.0245499}%
\special{ar 3304 3026 306 306  3.1423895 3.1816694}%
\special{ar 3304 3026 306 306  3.2995090 3.3387889}%
\special{ar 3304 3026 306 306  3.4566285 3.4959083}%
\special{ar 3304 3026 306 306  3.6137480 3.6530278}%
\special{ar 3304 3026 306 306  3.7708674 3.8101473}%
\special{ar 3304 3026 306 306  3.9279869 3.9672668}%
\special{ar 3304 3026 306 306  4.0851064 4.1243863}%
\special{ar 3304 3026 306 306  4.2422259 4.2815057}%
\special{ar 3304 3026 306 306  4.3993453 4.4386252}%
\special{ar 3304 3026 306 306  4.5564648 4.5957447}%
\special{ar 3304 3026 306 306  4.7135843 4.7528642}%
\special{ar 3304 3026 306 306  4.8707038 4.9099836}%
\special{ar 3304 3026 306 306  5.0278232 5.0671031}%
\special{ar 3304 3026 306 306  5.1849427 5.2242226}%
\special{ar 3304 3026 306 306  5.3420622 5.3813421}%
\special{ar 3304 3026 306 306  5.4991817 5.5384615}%
\special{ar 3304 3026 306 306  5.6563011 5.6955810}%
\special{ar 3304 3026 306 306  5.8134206 5.8527005}%
\special{ar 3304 3026 306 306  5.9705401 6.0098200}%
\special{ar 3304 3026 306 306  6.1276596 6.1669394}%
%
\special{pn 8}%
\special{ar 3638 2698 304 306  0.0000000 0.0393443}%
\special{ar 3638 2698 304 306  0.1573770 0.1967213}%
\special{ar 3638 2698 304 306  0.3147541 0.3540984}%
\special{ar 3638 2698 304 306  0.4721311 0.5114754}%
\special{ar 3638 2698 304 306  0.6295082 0.6688525}%
\special{ar 3638 2698 304 306  0.7868852 0.8262295}%
\special{ar 3638 2698 304 306  0.9442623 0.9836066}%
\special{ar 3638 2698 304 306  1.1016393 1.1409836}%
\special{ar 3638 2698 304 306  1.2590164 1.2983607}%
\special{ar 3638 2698 304 306  1.4163934 1.4557377}%
\special{ar 3638 2698 304 306  1.5737705 1.6131148}%
\special{ar 3638 2698 304 306  1.7311475 1.7704918}%
\special{ar 3638 2698 304 306  1.8885246 1.9278689}%
\special{ar 3638 2698 304 306  2.0459016 2.0852459}%
\special{ar 3638 2698 304 306  2.2032787 2.2426230}%
\special{ar 3638 2698 304 306  2.3606557 2.4000000}%
\special{ar 3638 2698 304 306  2.5180328 2.5573770}%
\special{ar 3638 2698 304 306  2.6754098 2.7147541}%
\special{ar 3638 2698 304 306  2.8327869 2.8721311}%
\special{ar 3638 2698 304 306  2.9901639 3.0295082}%
\special{ar 3638 2698 304 306  3.1475410 3.1868852}%
\special{ar 3638 2698 304 306  3.3049180 3.3442623}%
\special{ar 3638 2698 304 306  3.4622951 3.5016393}%
\special{ar 3638 2698 304 306  3.6196721 3.6590164}%
\special{ar 3638 2698 304 306  3.7770492 3.8163934}%
\special{ar 3638 2698 304 306  3.9344262 3.9737705}%
\special{ar 3638 2698 304 306  4.0918033 4.1311475}%
\special{ar 3638 2698 304 306  4.2491803 4.2885246}%
\special{ar 3638 2698 304 306  4.4065574 4.4459016}%
\special{ar 3638 2698 304 306  4.5639344 4.6032787}%
\special{ar 3638 2698 304 306  4.7213115 4.7606557}%
\special{ar 3638 2698 304 306  4.8786885 4.9180328}%
\special{ar 3638 2698 304 306  5.0360656 5.0754098}%
\special{ar 3638 2698 304 306  5.1934426 5.2327869}%
\special{ar 3638 2698 304 306  5.3508197 5.3901639}%
\special{ar 3638 2698 304 306  5.5081967 5.5475410}%
\special{ar 3638 2698 304 306  5.6655738 5.7049180}%
\special{ar 3638 2698 304 306  5.8229508 5.8622951}%
\special{ar 3638 2698 304 306  5.9803279 6.0196721}%
\special{ar 3638 2698 304 306  6.1377049 6.1770492}%
%
\special{pn 13}%
\special{pa 3324 3016}%
\special{pa 3608 3008}%
\special{fp}%
\special{sh 1}%
\special{pa 3608 3008}%
\special{pa 3540 2990}%
\special{pa 3554 3010}%
\special{pa 3542 3030}%
\special{pa 3608 3008}%
\special{fp}%
%
\special{pn 13}%
\special{pa 3648 2688}%
\special{pa 3932 2682}%
\special{fp}%
\special{sh 1}%
\special{pa 3932 2682}%
\special{pa 3866 2664}%
\special{pa 3880 2682}%
\special{pa 3866 2704}%
\special{pa 3932 2682}%
\special{fp}%
\put(19.0400,-19.0500){\makebox(0,0)[lt]{$\widetilde{\xi^j}$}}%
\put(40.7100,-19.3300){\makebox(0,0)[lt]{$\widetilde{\xi^j_{t_j}}$}}%
\put(15.4000,-16.3000){\makebox(0,0)[lt]{$e$}}%
\put(37.7600,-16.4000){\makebox(0,0)[lt]{$e$}}%
\put(24.0900,-30.5500){\makebox(0,0)[lt]{$e$}}%
%
\special{pn 8}%
\special{pa 3240 1030}%
\special{pa 3172 1246}%
\special{dt 0.045}%
\special{sh 1}%
\special{pa 3172 1246}%
\special{pa 3210 1190}%
\special{pa 3188 1196}%
\special{pa 3172 1176}%
\special{pa 3172 1246}%
\special{fp}%
\put(32.0000,-10.0000){\makebox(0,0)[lb]{$\Sigma$}}%
%
\special{pn 8}%
\special{pa 5430 1038}%
\special{pa 5360 1256}%
\special{dt 0.045}%
\special{sh 1}%
\special{pa 5360 1256}%
\special{pa 5400 1198}%
\special{pa 5376 1206}%
\special{pa 5362 1186}%
\special{pa 5360 1256}%
\special{fp}%
\put(54.0900,-10.0100){\makebox(0,0)[lb]{$\Sigma$}}%
\put(40.6500,-24.8300){\makebox(0,0)[lb]{$\Sigma$}}%
%
\special{pn 8}%
\special{pa 4126 2494}%
\special{pa 4074 2638}%
\special{dt 0.045}%
\special{sh 1}%
\special{pa 4074 2638}%
\special{pa 4116 2582}%
\special{pa 4092 2588}%
\special{pa 4078 2568}%
\special{pa 4074 2638}%
\special{fp}%
%
\special{pn 8}%
\special{pa 1522 1596}%
\special{pa 1522 1322}%
\special{dt 0.045}%
\special{sh 1}%
\special{pa 1522 1322}%
\special{pa 1502 1388}%
\special{pa 1522 1374}%
\special{pa 1542 1388}%
\special{pa 1522 1322}%
\special{fp}%
\put(15.0000,-12.8000){\makebox(0,0)[rb]{$H_{\lambda_j}$}}%
\put(27.9300,-33.8500){\makebox(0,0)[lt]{$E^j$}}%
%
\special{pn 20}%
\special{sh 1}%
\special{ar 2440 1600 10 10 0  6.28318530717959E+0000}%
\special{sh 1}%
\special{ar 2440 1600 10 10 0  6.28318530717959E+0000}%
%
\special{pn 20}%
\special{sh 1}%
\special{ar 4640 1610 10 10 0  6.28318530717959E+0000}%
\special{sh 1}%
\special{ar 4640 1610 10 10 0  6.28318530717959E+0000}%
\put(24.6000,-9.5000){\makebox(0,0)[rb]{$\gamma_{\xi^j}(t_j)$}}%
%
\special{pn 8}%
\special{pa 2290 990}%
\special{pa 2430 1590}%
\special{dt 0.045}%
\special{sh 1}%
\special{pa 2430 1590}%
\special{pa 2434 1522}%
\special{pa 2418 1538}%
\special{pa 2396 1530}%
\special{pa 2430 1590}%
\special{fp}%
\put(46.8000,-9.8000){\makebox(0,0)[rb]{$\gamma_{\xi^j}(t_j)$}}%
%
\special{pn 8}%
\special{pa 4500 1000}%
\special{pa 4640 1600}%
\special{dt 0.045}%
\special{sh 1}%
\special{pa 4640 1600}%
\special{pa 4644 1532}%
\special{pa 4628 1548}%
\special{pa 4606 1540}%
\special{pa 4640 1600}%
\special{fp}%
\put(33.4000,-23.8000){\makebox(0,0)[rb]{$\gamma_{\xi^j}(t_j)$}}%
%
\special{pn 8}%
\special{pa 3170 2420}%
\special{pa 3310 3020}%
\special{dt 0.045}%
\special{sh 1}%
\special{pa 3310 3020}%
\special{pa 3314 2952}%
\special{pa 3298 2968}%
\special{pa 3276 2960}%
\special{pa 3310 3020}%
\special{fp}%
%
\special{pn 20}%
\special{sh 1}%
\special{ar 3310 3020 10 10 0  6.28318530717959E+0000}%
\special{sh 1}%
\special{ar 3310 3020 10 10 0  6.28318530717959E+0000}%
\put(36.8000,-13.1000){\makebox(0,0)[rb]{$H_{\lambda_j}$}}%
%
\special{pn 8}%
\special{pa 3740 1610}%
\special{pa 3740 1320}%
\special{dt 0.045}%
\special{sh 1}%
\special{pa 3740 1320}%
\special{pa 3720 1388}%
\special{pa 3740 1374}%
\special{pa 3760 1388}%
\special{pa 3740 1320}%
\special{fp}%
\put(18.6000,-16.3000){\makebox(0,0)[lt]{$\xi^j$}}%
\put(40.6000,-16.3000){\makebox(0,0)[lt]{$\xi^j$}}%
\put(27.0000,-30.5000){\makebox(0,0)[lt]{$\xi^j$}}%
%
\special{pn 8}%
\special{pa 3750 1610}%
\special{pa 4030 1610}%
\special{dt 0.045}%
\special{sh 1}%
\special{pa 4030 1610}%
\special{pa 3964 1590}%
\special{pa 3978 1610}%
\special{pa 3964 1630}%
\special{pa 4030 1610}%
\special{fp}%
%
\special{pn 8}%
\special{pa 2370 3020}%
\special{pa 2370 2730}%
\special{dt 0.045}%
\special{sh 1}%
\special{pa 2370 2730}%
\special{pa 2350 2798}%
\special{pa 2370 2784}%
\special{pa 2390 2798}%
\special{pa 2370 2730}%
\special{fp}%
\put(23.2000,-27.3000){\makebox(0,0)[rb]{$H_{\lambda_j}$}}%
\end{picture}%
\hspace{0.5truecm}}

\vspace{0.5truecm}

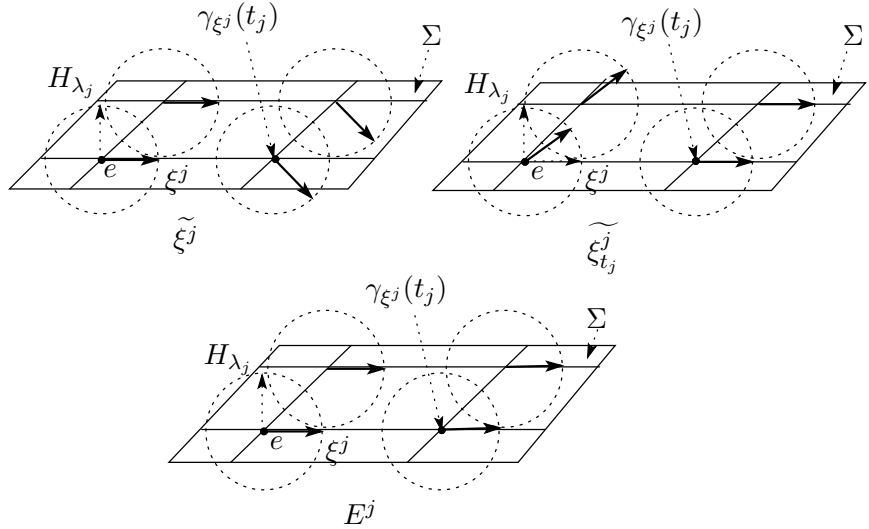
\centerline{{\bf Figure 4.}}

\vspace{0.5truecm}

By using these lemmas, we prove Theorem A.  

\vspace{0.5truecm}

\noindent
{\it Proof of Theorem A.} In this proof, we use the notations as in Example 2.  
Set $M_{x_{\xi_0,t_1,\cdots,t_k}}:=S_{\mathfrak b,\overline{\it l}_1,\cdots,\overline{\it l}_k}
\cdot x_{\xi_0,t_1,\cdots,t_k}$.  
Denote by $H^{x_{\xi_0,t_1,\cdots,t_k}}$ the mean curvature vector field of $M_{x_{\xi_0,t_1,\cdots,t_k}}$.  
Let $\{e^0_1,\cdots,e^0_{m_0}\}$ be an orthonormal base of $\mathfrak b$ and 
$(H_{\lambda})_{\mathfrak b}=\sum_{i=1}^{m_0}H_{\lambda}^ie_i^0$ be the $\mathfrak b$-component of 
$H_{\lambda}$.  
According to the fact (iv) stated in Introduction, we have 
$$M_{x_{\xi_0,t_1,\cdots,t_k}}=L_{x_{\xi_0,t_1,\cdots,t_k}}(S_{\mathfrak b,{\it l}_{\xi^1,t_1},\cdots,
{\it l}_{\xi^k,t_k}}\cdot e).$$
Denote by $\widehat H^{\xi_0,t_1,\cdots,t_k}$ the mean curvature vector field of 
$S_{\mathfrak b,{\it l}_{\xi^1,t_1},\cdots,{\it l}_{\xi^k,t_k}}\cdot e$.  
According to Lemma 3.1, we have 
$$(\widehat H^{\xi_0,t_1,\cdots,t_k})_e=\sum_{\lambda\in\triangle_+}m_{\lambda}(H_{\lambda})_{\mathfrak b}
-\sum_{i=1}^k\vert\vert\lambda_i\vert\vert\tanh(\vert\vert\lambda_i\vert\vert t_i)
(m_{\lambda_i}+2m_{2\lambda_i})\xi^i_{t_i}$$
and hence 
$$\begin{array}{l}
\displaystyle{(H^{x_{\xi_0,t_1,\cdots,t_k}})_{x_{\xi_0,t_1,\cdots,t_k}}
=\sum_{\lambda\in\triangle_+}\sum_{i=1}^{m_0}m_{\lambda}H_{\lambda}^i(E^0_i)_{x_{\xi_0,t_1,\cdots,t_k}}}\\
\hspace{4.2truecm}\displaystyle{-\sum_{i=1}^k
\vert\vert\lambda_i\vert\vert\tanh(\vert\vert\lambda_i\vert\vert t_i)
(m_{\lambda_i}+2m_{2\lambda_i})(E^i)_{x_{\xi_0,t_1,\cdots,t_k}}.}
\end{array}
\leqno{(4.7)}$$
Define a tangent vector field $Z$ over $\Sigma$ by $Z_x:=(H^x)_x\,\,(x\in\Sigma)$.  
According to $(4.7)$, we have 
$$\begin{array}{l}
\displaystyle{Z_{x_{\xi_0,t_1,\cdots,t_k}}
=\sum_{\lambda\in\triangle_+}\sum_{i=1}^{m_0}m_{\lambda}H_{\lambda}^i(E^0_i)_{x_{\xi_0,t_1,\cdots,t_k}}}\\
\hspace{2.5truecm}\displaystyle{-\sum_{i=1}^k\vert\vert\lambda_i\vert\vert\tanh(\vert\vert\lambda_i\vert\vert t_i)
(m_{\lambda_i}+2m_{2\lambda_i})(E^i)_{x_{\xi_0,t_1,\cdots,t_k}}.}
\end{array}\leqno{(4.8)}$$
Define a coordinate 
$\phi=(u_1,\cdots,u_{m_0+k}):\Sigma\to{\mathbb R}^{m_0+k}$ of $\Sigma$ by 
$$\phi(x_{\sum_{i=1}^{m_0}s_ie^0_i,t_1,\cdots,t_k}):=(s_1,\cdots,s_{m_0},t_1,\cdots,t_k)$$
($s_1,\cdots,s_{m_0},t_1,\cdots,t_k\in\mathbb R$).  
We can show $\frac{\partial}{\partial u_i}=E^0_i\,\,(i=1,\cdots,m_0)$ and 
$\frac{\partial}{\partial u_{m_0+j}}=E^j\,\,(j=1,\cdots,k)$.  
Hence $\phi$ is a Euclidean coordinate of $\Sigma$.  
Under the identification of $\Sigma$ and ${\mathbb R}^{m_0+k}$ by $\phi$, we regard $Z$ as a tangent vector field 
on ${\mathbb R}^{m_0+k}$.  
Then $Z$ is described as 
$$\begin{array}{l}
\displaystyle{Z_{(u_1,\cdots,u_{m_0+k})}
=(\sum_{\lambda\in\triangle_+}m_{\lambda}H_{\lambda}^1,\cdots,
\sum_{\lambda\in\triangle_+}m_{\lambda}H_{\lambda}^{m_0},}\\
\hspace{3truecm}\displaystyle{-\vert\vert\lambda_1\vert\vert\tanh(\vert\vert\lambda_1\vert\vert u_{m_0+1})
(m_{\lambda_1}+2m_{2\lambda_1}),}\\
\hspace{3truecm}\displaystyle{\cdots,-\vert\vert\lambda_k\vert\vert\tanh(\vert\vert\lambda_k\vert\vert u_{m_0+k})
(m_{\lambda_k}+2m_{2\lambda_k})).}
\end{array}\leqno{(4.9)}$$
Fix $(a_1,\cdots,a_{m_0},t_1,\cdots,t_k)\in{\mathbb R}^{m_0+k}$.  
Let $c$ be the integral curve of $Z$ starting from 
$(a_1,\cdots,a_{m_0},t_1,\cdots,t_k)$ and let $c=(c_1,\cdots,c_{m_0+k})$.  
We suffice to investigate $c$ to investigate the mean curvature flow 
starting from $M_{x_{\sum_{i=1}^{m_0}a_ie^0_i,t_1,\cdots,t_k}}$
From $c'(t)=Z_{c(t)}$, we have 
$c'_i(t)=\sum_{\lambda\in\triangle_+}m_{\lambda}H_{\lambda}^i$ ($i=1,\cdots,m_0$) and 
$c'_{m_0+j}(t)=-(m_{\lambda_j}+2m_{2\lambda_j})\vert\vert\lambda_j\vert\vert\tanh$\newline
$(\vert\vert\lambda_j\vert\vert c_{m_0+j}(t))$ ($j=1,\cdots,k$).  
By solving $c'_i(t)=\sum_{\lambda\in\triangle_+}m_{\lambda}H_{\lambda}^i$ under the initial condition 
$c_i(0)=a_i$, we have 
$$c_i(t)=a_i+t\sum_{\lambda\in\triangle_+}m_{\lambda}H_{\lambda}^i.\leqno{(4.10)}$$
Also, by solving $c'_{m_0+j}(t)=-(m_{\lambda_j}+2m_{2\lambda_j})\vert\vert\lambda_j\vert\vert
\tanh(\vert\vert\lambda_j\vert\vert c_{m_0+j}(t))$ under the initial condition $c_{m_0+j}(0)=t_j$, we have 
$$c_{m_0+j}(t)=\frac{1}{\vert\vert\lambda_j\vert\vert}
{\rm arcsinh}\left(e^{-\vert\vert\lambda_j\vert\vert^2(m_{\lambda_j}+2m_{2\lambda_j})t}
\sinh(\vert\vert\lambda_j\vert\vert t_j)\right).\leqno{(4.11)}$$
From $(4.10)$ and $(4.11)$, we can derive $T=\infty$, $\lim\limits_{t\to\infty}\sum_{i=1}^{m_0}c_i(t)^2
=\infty$ ($i=1,\cdots,m_0$) and $\lim\limits_{t\to\infty}c_{m_0+j}(t)=0$ ($j=1,\cdots,k$).  
If $t_1=\cdots=t_k=0$, then we have $c_{m_0+j}\equiv 0$ ($j=1,\cdots,m_0$).  
Hence the mean curvature flow starting from $M_{x_{\xi_0,0,\cdots,0}}$ 
($x_{\xi_0,0,\cdots,0}\in{\rm Exp}(\mathfrak b)$) consists of the leaves of 
${\mathfrak F}_{\mathfrak b,\overline{\it l}_1,\cdots,\overline{\it l}_k}$ through points of 
${\rm Exp}(\mathfrak b)$.  Also, according to the fact (iv) stated in Introduction, 
the leaves of ${\mathfrak F}_{\mathfrak b,\overline{\it l}_1,\cdots,\overline{\it l}_k}$ 
through points of ${\rm Exp}(\mathfrak b)$ are congruent to 
$S_{\mathfrak b,\overline{\it l}_1,\cdots,\overline{\it l}_k}\cdot e$.  
Therefore, the mean curvature flow starting from $M_{x_{\xi_0,0,\cdots,0}}$ is self-similar.  
From $\lim\limits_{t\to\infty}\sum_{i=1}^{m_0}c_i(t)^2=\infty$ ($i=1,\cdots,m_0$) and 
$\lim\limits_{t\to\infty}c_{m_0+j}(t)=0$ ($j=1,\cdots,k$), we see that 
the mean curvature flow starting from any leaf of 
${\mathfrak F}_{\mathfrak b,\overline{\it l}_1,\cdots,\overline{\it l}_k}$ asymptotes the mean curvature flow 
starting from the leaf of ${\mathfrak F}_{\mathfrak b,\overline{\it l}_1,\cdots,\overline{\it l}_k}$ passing through 
a point of ${\rm Exp}(\mathfrak b)$.  
\hspace{1.5truecm}q.e.d.

\vspace{0.5truecm}

According to this proof, we obtain the following fact.  

\vspace{0.5truecm}

\noindent
{\bf Corollary 4.1.} {\sl {\rm(i)} The mean curvature flow starting from $M_{x_{\xi_0,0,\cdots,0}}$ is 
self-similar.  

{\rm(ii)} The mean curvature flow starting from $M_{x_{\xi_0,t_1,\cdots,t_k}}$ 
($(t_1,\cdots,t_k)\not=(0,\cdots,0)$) asymptotes the flow starting from 
$M_{x_{\xi_0,0,\cdots,0}}$.  In more detail, the distance between $M_{x_{\xi_0,t_1,\cdots,t_k}}$ and 
$M_{x_{\xi_0,0,\cdots,0}}$ is equal to 
$$\sqrt{\sum_{j=1}^k\frac{1}{\vert\vert\lambda_j\vert\vert^2}
{\rm arcsinh}^2\left(e^{-\vert\vert\lambda_j\vert\vert^2(m_{\lambda_j}+2m_{2\lambda_j})t}
\sinh(\vert\vert\lambda_j\vert\vert t_j)\right),}$$
which converges to zero as $t\to\infty$.}

\vspace{0.5truecm}

Next we prove Theorem B.  

\vspace{0.5truecm}

\noindent
{\it Proof of Theorem B.} In case of $\mathfrak b=\{0\}$, the relation $(4.9)$ is as follows:
$$\begin{array}{l}
\displaystyle{Z_{(u_1,\cdots,u_k)}
=(-\vert\vert\lambda_1\vert\vert\tanh(\vert\vert\lambda_1\vert\vert u_{m_0+1})
(m_{\lambda_1}+2m_{2\lambda_1}),}\\
\hspace{3truecm}\displaystyle{\cdots,-\vert\vert\lambda_k\vert\vert\tanh(\vert\vert\lambda_k\vert\vert u_{m_0+k})
(m_{\lambda_k}+2m_{2\lambda_k})).}
\end{array}\leqno{(4.12)}$$
Hence, according to the dicussion in the proof of Theorem A, 
the mean curvature flow starting from any leaf of 
${\mathfrak F}_{\mathfrak b,\overline{\it l}_1,\cdots,\overline{\it l}_k}$ converges to the only minimal leaf 
$S_{\mathfrak b,\overline{\it l},\cdots,\overline{\it l}_k}\cdot e$.  
Furthermore, the flow converges to the minimal leaf in $C^{\infty}$-topology because the flow consists of 
$S_{\mathfrak b,\overline{\it l}_1,\cdots,\overline{\it l}_k}$-orbits and the limit submanifold also is 
a $S_{\mathfrak b,\overline{\it l}_1,\cdots,\overline{\it l}_k}$-orbit.  \hspace{2.4truecm}q.e.d.

\vspace{0.8truecm}

\centerline{{\bf References}}

\vspace{0.5truecm}

{\small
\noindent
[AB] M. M. Alexandrino and M. Radeschi, Mean curvature flow of singular Riemannian foliations, 

J. Geom. Anal. (to appear).

\noindent
[AB] B. Andrews and C. Baker, Mean curvature flow of pinched submanifolds to spheres, J. Diffe-

rential Geom. {\bf 85} (2010) 357-396.

\noindent
[H1] G. Huisken, Flow by mean curvature of convex surfaces into spheres, J. Differential Geom. 

{\bf 20} (1984) 237-266.

\noindent
[H2] G. Huisken, Contracting convex hypersurfaces in Riemannian manifolds by their mean cur-

vature, Invent. math. {\bf 84} (1986) 463-480.

\noindent
[K1] N. Koike, Examples of a complex hyperpolar action without singular orbit, Cubo A Math. 

J. {\bf 12} (2010) 131-147.

\newpage

\noindent
[K2] N. Koike, Collapse of the mean curvature flow for equifocal submanifolds, Asian J. Math. 

{\bf 15} (2011) 101-128.

\noindent
[K3] N. Koike, Collapse of the mean curvature flow for isoparametric submanifolds in a symme-

tric space of non-compact type, Kodai Math. J. {\bf 37} (2014) 355-382.

\noindent
[M] J. Milnor, 
Curvatures of left invariant metrics on Lie groups, Adv. Math. {\bf 21} 
(1976) 293--329.
}

\vspace{0.5truecm}

\rightline{Department of Mathematics, Faculty of Science, }
\rightline{Tokyo University of Science}
\rightline{1-3 Kagurazaka Shinjuku-ku,}
\rightline{Tokyo 162-8601, Japan}
\rightline{(e-mail: koike@ma.kagu.tus.ac.jp)}

\end{document}